\documentclass[11pt]{amsart}
\usepackage{latexsym}
\usepackage{amsmath}
\usepackage{amssymb}

\newcommand{\eproof}{\mbox{\ }\hfill $\Box$ \par \vskip 10pt}

\newtheorem{Theorem}{Theorem}[section]
\newtheorem{lemma}[Theorem]{Lemma}
\newtheorem{prop}[Theorem]{Proposition}

\numberwithin{equation}{section}

\def\cal{\mathcal}

\begin{document}

\title[Improved parametrix in the glancing region]{Improved parametrix in the glancing region for the interior Dirichlet-to-Neumann map}

\author[G. Vodev]{Georgi Vodev}

\address {Universit\'e de Nantes, Laboratoire de Math\'ematiques Jean Leray, 2 rue de la Houssini\`ere, BP 92208, 44322 Nantes Cedex 03, France}
\email{Georgi.Vodev@univ-nantes.fr}

\date{}

\begin{abstract} We study the semi-classical microlocal structure of the Dirichlet-to-Neumann map for an arbitrary 
compact Riemannian manifold with a non-empty smooth boundary. We build a new, improved parametrix in the glancing region compaired with that
one built in \cite{kn:V1}, \cite{kn:V4}. We also study the way in which the parametrix depends on the refraction index. As a consequence, 
 we improve the transmission eigenvalue-free regions obtained in \cite{kn:V4} in the isotropic case when the restrictions of the refraction
indices on the boundary coincide. 
\end{abstract} 

\maketitle

\setcounter{section}{0}
\section{Introduction and statement of results}

Let $(X,{\cal G})$ be a compact Riemannian manifold of dimension $d={\rm dim}\, X\ge 2$ with a non-empty smooth boundary $\partial X$
and let $\Delta_X$ denote the negative Laplace-Beltrami operator on
$(X,{\cal G})$.  
Given a function $f\in H^{m+1}(\partial X)$, let  $u$ solve the equation
\begin{equation}\label{eq:1.1}
\left\{
\begin{array}{lll}
 \left(\Delta_X+\lambda^2n(x)\right)u=0&\mbox{in}& X,\\
 u=f&\mbox{on}&\partial X,
\end{array}
\right.
\end{equation}
where $\lambda\in{\bf C}$, $1\ll|{\rm Im}\,\lambda|\ll {\rm Re}\,\lambda$ and $n\in C^\infty(\overline X)$ is a strictly positive function
called refraction index. 
The Dirichlet-to-Neumann (DN) map 
$${\cal N}(\lambda;n):H^{m+1}(\partial X)\to H^m(\partial X)$$
is defined by
$${\cal N}(\lambda;n)f:=\partial_\nu u|_{\partial X}$$
 where $\nu$ is the unit inner normal to $\partial X$. Introduce the semi-classical parameter $0<h\ll 1$ such that 
 ${\rm Re}(h\lambda)^2=1$ and set $z=(h\lambda)^2=1+i{\rm Im}\,z$ with $0<|{\rm Im}\,z|\le 1$. Then the problem (\ref{eq:1.1}) can be 
 rewritten as follows
 \begin{equation}\label{eq:1.2}
\left\{
\begin{array}{lll}
 \left(h^2\Delta_X+zn(x)\right)u=0&\mbox{in}& X,\\
 u=f&\mbox{on}&\partial X.
\end{array}
\right.
\end{equation}
 Define the semi-classical DN map, $N(h,z)$, by
 $$N(h,z)f:={\cal D}_\nu u|_{\partial X}=-ih{\cal N}(\lambda;n)f$$
 where ${\cal D}_\nu:=-ih\partial_\nu$. 
 Denote by $\Delta_{\partial X}$ the negative Laplace-Beltrami operator on $(\partial X,{\cal G}_0)$, which is
  a Riemannian manifold without boundary of dimension $d-1$, where ${\cal G}_0$ is the Riemannian
metric on $\partial X$ induced by the metric ${\cal G}$. Let $r_0(x',\xi')\ge 0$ be the principal symbol of 
$-\Delta_{\partial X}$ written in the coordinates $(x',\xi')\in T^*\partial X$. The glancing region, $\Sigma$, for the problem (\ref{eq:1.2}) 
(resp. (\ref{eq:1.1})) is defined by
$$\Sigma:=\{(x',\xi')\in T^*\partial X:r_0(x',\xi')=n_0(x')\},\quad n_0:=n|_{\partial X}.$$
Our goal in the present paper is to build a semi-classical parametrix for the operator $N(h,z)$ in a neighbourhood of $\Sigma$
for $|{\rm Im}\,z|\ge h^{2/3-\epsilon}$, $0<\epsilon\ll 1$ being arbitrary. Since $h\sim|\lambda|^{-1}$, it is easy to see that on the
$\lambda-$ plane this region takes the form $|{\rm Im}\,\lambda|\ge |\lambda|^{1/3+\epsilon}$, $|\lambda|\gg 1$. 
Note that such a parametrix has been previously constructed
in \cite{kn:V1}, \cite{kn:V4} for $|{\rm Im}\,z|\ge h^{1/2-\epsilon}$ (this corresponds to the region $|{\rm Im}\,\lambda|\ge |\lambda|^{1/2+\epsilon}$, $|\lambda|\gg 1$, on the $\lambda-$ plane). Roughly speaking, the smaller $|{\rm Im}\,z|$ is, the harder is to
construct a parametrix for $N(h,z)$. Note also that a semi-classical parametrix for the operator $N(h,z)$ in a neighbourhood of $\Sigma$
has been built in \cite{kn:V2} for $|{\rm Im}\,z|\ge h^{1-\epsilon}$ (which corresponds to the region $|{\rm Im}\,\lambda|\ge |\lambda|^{\epsilon}$, $|\lambda|\gg 1$, on the $\lambda-$ plane) but under the additional assumption that the boundary $\partial X$
is strictly concave. Under this condition another semi-classical parametrix was built in \cite{kn:S} for $|{\rm Im}\,z|\sim h^{2/3}$.

It has been shown in \cite{kn:V1} that $N(h,z)\in OPS^{0,1}_{\delta,\delta}(\partial X)$ (see Section 2 for the definition of the $h-\Psi$DOs
of class $OPS^{k_1,k_2}_{\delta_1,\delta_2}$) for $|{\rm Im}\,z|\ge h^{1/2-\epsilon}$ with
$\delta=1/2-\epsilon$, and a principal symbol, $\rho$, defined by
$$\rho(x',\xi';z)=\sqrt{-r_0(x',\xi')+zn_0(x')},\quad{\rm Im}\,\rho>0.$$
Moreover, outside the glancing region the operator $N(h,z)$ belongs to a much better class, due to the fact that there 
$|\rho|$ is lower bounded by a 
positive constant. To be more precise, we choose a cut-off function $\chi\in C_0^\infty(T^*\partial X)$, independent of $h$ and $z$, such that
$\chi=1$ in a small, $h-$ independent neighbourhood of $\Sigma$. It follows from the analysis in \cite{kn:V1} that $N(h,z){\rm Op}_h(1-\chi)\in OPS^{0,1}_{0,0}(\partial X)$ for $|{\rm Im}\,z|\ge h^{1-\epsilon}$ with a principal symbol $\rho(1-\chi)$. In other words, the
condition $|{\rm Im}\,z|\ge h^{1/2-\epsilon}$ is only required in \cite{kn:V1} to study the operator $N(h,z)$ near $\Sigma$. Note that the
full symbol of $N(h,z)$ depends on the functions $n_k=\partial_\nu^k n|_{\partial X}$, $k=0,1,...,$ and their derivatives. The way in which 
the parametrix of the operator $N(h,z)$ depends on $n_k$ is studied in \cite{kn:V4}. Near the glancing region the analysis in \cite{kn:V4} again requires
the condition $|{\rm Im}\,z|\ge h^{1/2-\epsilon}$, while outside $\Sigma$ it works for $|{\rm Im}\,z|\ge h^{1-\epsilon}$. 

In the present paper we will extend the analysis near $\Sigma$ to the larger region $|{\rm Im}\,z|\ge h^{2/3-\epsilon}$. 
The first difficulty to deal with is to give a reasonable definition of the operator ${\rm Op}_h(\rho)$ when 
$|{\rm Im}\,z|\sim h^{2/3-\epsilon}$. Indeed, in this case $\chi\rho\in S^{0,0}_{\delta,\delta}(\partial X)$ with $\delta=2/3-\epsilon>1/2$
and we do not have a good calculus for $h-\Psi$DOs with symbols in this class. To overcome this difficulty we use the 
second microlocalization with respect to $\Sigma$. Note that this approach proved very usefull when studying the 
resonances near cubic curves in the case of scattering by strictly convex obstacles (see \cite{kn:SZ}).
It has been also successfully used in \cite{kn:G} to study the location of the resonances for various exterior transmission problems
associated to transparent strictly convex obstacles. It has already been used in the contex of the interior DN map in \cite{kn:S}
and \cite{kn:V2} to build a parametrix when the boundary $\partial X$ is strictly concave. Here we will use it for an
arbitrary Riemannian manifold. Roughly speaking, it consists of using $h-$FIOs to transform our
problem (\ref{eq:1.2}) microlocally
near the boundary into a simpler equation (see the model equation (\ref{eq:2.3}) in Section 2) for which it is 
easier to construct a microlocal
parametrix. Then the global parametrix is obtained by using a suitable partition of the unity on $\Sigma$. More precisely,
let ${\cal W}\subset T^*\partial X$ be a small neighbourhood of $\Sigma$ such that supp$\,\chi\subset {\cal W}$
 and  cover ${\cal W}$ with a finite number of sufficiently small, open domains, ${\cal W}_j\subset T^*\partial X$, $j=1,...,J$, 
$\Sigma\subset{\cal W}\subset\cup_{j=1}^J{\cal W}_j$. Choose functions $\chi_j,\psi_j\in C_0^\infty({\cal W}_j)$ such that $\psi_j=1$
on supp$\,\chi_j$ and $\chi=\sum_{j=1}^J\chi_j$. It is well-known (e.g. see Theorem 12.3 of \cite{kn:Z}) that there are an open, bounded domain
$Y_j\subset{\bf R}^{d-1}$ and a symplectomorphism $\kappa_j:{\cal W}_j\to T^*Y_j$ such that in the new coordinates
$(y,\eta)=\kappa_j(x',\xi')$ we have
$$\eta_1=n_0^{-1}r_0(x',\xi')-1$$
and $\kappa_j(\Sigma\cap{\cal W}_j)=\{\eta_1=0\}$. We can write the function $\rho$ as follows
$$\rho=n_0^{1/2}\varrho\circ\kappa_j^{-1}$$
where
$$\varrho(\eta_1)=\sqrt{-\eta_1+i{\rm Im}\,z},\quad {\rm Im}\,\varrho>0.$$
Let $U_j:L^2(\pi{\cal W_j})\to L^2(Y_j)$ be an elliptic, zero-order $h-$FIO associated to
$\kappa_j$, where $\pi:T^*\partial X\to \partial X$ denotes the projection $\pi(x',\xi')=x'$. Then we define the operators 
$$\widetilde{{\rm Op}}_h(\rho\chi)=\sum_{j=1}^J{\rm Op}_h(n_0^{1/2}\chi_j)U_j^{-1}{\rm Op}_h(\phi_j\varrho)U_j
{\rm Op}_h(\psi_j):L^2(\partial X)\to L^2(\partial X)$$
where $\phi_j(\eta_1)\in C_0^\infty$ is such that $\phi_j=1$ on supp$\,\chi_j\circ\kappa_j$, 
and 
$$\widetilde{{\rm Op}}_h(\rho)={\rm Op}_h(\rho(1-\chi))+\widetilde{{\rm Op}}_h(\rho\chi).$$
We can define similarly the operators $\widetilde{{\rm Op}}_h(\rho^k\chi)$ and $\widetilde{{\rm Op}}_h(\rho^k)$, 
$k$ being an arbitrary integer. When $|{\rm Im}\,z|\ge h^{1/2-\epsilon}$ one can see that the operator 
$\widetilde{{\rm Op}}_h(\rho\chi)$ coincides with the standard $h-\Psi$DO ${\rm Op}_h(\rho\chi)$ mod ${\cal O}(h^\epsilon)$.
We refer to Section 4
of \cite{kn:SZ}, Section 4 of \cite{kn:G} as well as Sections 10, 11 and 12 of \cite{kn:Z} for more information about the properties of these
operators. In the present paper we only make use of some very basic properties of the $h-\Psi$DOs and $h-$FIOs.

Thus we reduce the parametrix construction in the glancing region to building a parametrix for a model equation (see Section 3).
On the other hand, the parametrix construction for the model equation is carried out in Section 2 following that one in Section 6 of \cite{kn:V2}
in the case of strictly concave boundary. Note that the model equation in \cite{kn:V2} is much simpler than this one we study in the present paper.
This is due to the fact that the strict concavity condition allows to use the symplectic normal form proved in \cite{kn:PV}. No such normal
forms exist, however, in the general case and therefore we have to work with a model equation which is relatively complicated. 
Nevertheless, we show that the parametrix construction still works. 
As a consequence we get the following

\begin{Theorem} Let $|{\rm Im}\,z|\ge h^{2/3-\epsilon}$, $0<\epsilon\ll 1$. Then we have
\begin{equation}\label{eq:1.3}
\left\|N(h,z){\rm Op}_h(\chi)-\widetilde{{\rm Op}}_h(\rho\chi)\right\|
\lesssim h|{\rm Im}\,z|^{-1}.
\end{equation}
\end{Theorem}

Hereafter $\|\cdot\|$ denotes the $L^2(\partial X)\to L^2(\partial X)$ norm. 
This theorem is a significant improvement upon the results in \cite{kn:V1} and \cite{kn:V4}. Indeed, 
an analog of this theorem has been proved in \cite{kn:V1} and \cite{kn:V4} but for $|{\rm Im}\,z|\ge h^{1/2-\epsilon}$ 
and with an worse bound in the right-hand side of (\ref{eq:1.3}) (with $|{\rm Im}\,z|^{-3/2}$ in place of $|{\rm Im}\,z|^{-1}$). 
In fact, we get a full expansion in powers of $h$ of the operator 
$N{\rm Op}_h(\chi)$ and we study the way in which it depends on the functions $n_k$ (see Theorem 3.1).

Outside the glancing region we have a better bound for the DN map. Indeed, it has been proved in
\cite{kn:V1} and \cite{kn:V4} that for $|{\rm Im}\,z|\ge h^{1-\epsilon}$ the following estimate holds:
\begin{equation}\label{eq:1.4}
\left\|N(h,z){\rm Op}_h(1-\chi)-{\rm Op}_h(\rho(1-\chi))\right\|
\lesssim h.
\end{equation}
 Combining (\ref{eq:1.3}) and (\ref{eq:1.4})
we get the following

\begin{Theorem} Let $|{\rm Im}\,z|\ge h^{2/3-\epsilon}$, $0<\epsilon\ll 1$. Then we have
\begin{equation}\label{eq:1.5}
\left\|N(h,z)-\widetilde{{\rm Op}}_h(\rho)\right\|
\lesssim h|{\rm Im}\,z|^{-1}.
\end{equation}
\end{Theorem}

This theorem provides a good approximation of the DN map by an $h-\Psi$DO as long as $|{\rm Im}\,z|\ge h^{2/3-\epsilon}$.
Even better approximations are given in Theorem 3.5. 
For many applications, however, one needs to have some less accurate approximation of the DN map but for smaller $|{\rm Im}\,z|$.
Indeed, such an approximation has been proved in \cite{kn:V3} for $|{\rm Im}\,z|\ge Ch$, provided the constant $C>0$ is taken big enough.
Having high-frequency approximations of the DN map proves very usefull when studying the location of the complex eigenvalues 
associated to boundary value problems with dissipative boundary conditions or to interior transmission problems. In particular,
this proves crucial to get parabolic transmission eigenvalue-free regions (see \cite{kn:V1}, \cite{kn:V2}, \cite{kn:V3}, \cite{kn:V4}). As an application
of our parametrix we improve the transmission eigenvalue-free region obtained in \cite{kn:V4} in the case of the degenerate
isotropic interior transmission problem (see Theorem 4.1).

\section{Parametrix construction for the model equation} 

Let $Y\subset{\bf R}^{d-1}$, $d\ge 2$, be a bounded, open domain. Given $k_1,k_2\in{\bf R}$, $\delta_1,\delta_2\ge 0$, define the class of symbols
$S^{k_1,k_2}_{\delta_1,\delta_2}(Y)$ as being the set of all functions $a\in C^\infty(T^*Y)$, supp$_y\,a(y,\eta)\subset Y$, depending on a semi-classical parameter
$0<h\ll 1$ and satisfying the bounds
$$
\left|\partial_{y}^\alpha\partial_{\eta}^\beta a(y,\eta)\right|\le C_{\alpha,\beta}h^{-k_1-\delta_1|\alpha|-\delta_2|\beta|}
\langle\eta\rangle^{k_2-|\beta|}
$$
for all multi-indices $\alpha$ and $\beta$. We then define the $h-\Psi$DO with a symbol $a$ by
$$a(y,{\cal D}_y)f:=\left({\rm Op}_h(a)f\right)(y):=(2\pi h)^{-d+1}\int_{T^*Y}e^{-\frac{i}{h}\langle y-w,\eta\rangle}a(y,\eta)f(w)dwd\eta$$
where ${\cal D}_y:=-ih\partial_y$. We will denote by $OPS^{k_1,k_2}_{\delta_1,\delta_2}(Y)$ the set of all operators ${\rm Op}_h(a)$
with symbols $a\in S^{k_1,k_2}_{\delta_1,\delta_2}(Y)$. 
Since $Y$ is bounded, we have the following simple criteria for an $h-\Psi$DO to be bounded on $L^2(Y)$:
\begin{equation}\label{eq:2.1}
\left\|{\rm Op}_h(a)\right\|_{L^2(Y)\to L^2(Y)}\le C\sum_{0\le|\alpha|\le d}\sup_{y,\eta}|\partial_y^\alpha a(y,\eta)|
\end{equation}
where $C>0$ is a constant independent of $a$ and $h$. In particular, (\ref{eq:2.1}) implies that if $a\in S^{k_1,0}_{0,\delta_2}(Y)$, then
${\rm Op}_h(a)={\cal O}(h^{-k_1}):{L^2(Y)\to L^2(Y)}$. Note that (\ref{eq:2.1}) is no longer true if $Y={\bf R}^{d-1}$. Indeed, in that case one
also needs some information concerning the derivatives with respect to the variable $\eta$.

Given two symbols $a$ and $b$ and an integer $M\ge 1$, set
$$E_M(a,b)=\sum_{0\le|\alpha|\le M}\frac{(-ih)^{|\alpha|}}{|\alpha|!}\partial_\eta^\alpha a(y,\eta)\partial_y^\alpha b(y,\eta).$$
It is easy to see that if $a\in S^{k_1,k_2}_{\delta_1,\delta_2}(Y)$, $b\in S^{k'_1,k'_2}_{\delta'_1,\delta'_2}(Y)$ 
with $\delta_2+\delta'_1\le 1$, then 
 $E_M(a,b)\in S^{k^\sharp_1,k^\sharp_2}_{\delta^\sharp_1,\delta^\sharp_2}(Y)$, where $k^\sharp_j=k_j+k'_j$,
$\delta^\sharp_j=\max\{\delta_j,\delta'_j\}$, $j=1,2$.
The following proposition follows from the calculus developed in Section 7 of \cite{kn:DS}.

\begin{prop} Let $a\in S^{k_1,0}_{\delta_1,\delta_2}(Y)$ and let $b\in C^\infty(T^*Y)$ satisfy
$$
\left|\partial_{y}^\alpha b(y,\eta)\right|\le C_{\alpha}h^{-k_2-\delta'_1|\alpha|}
$$
for all multi-indices $\alpha$. If $\delta_2+\delta'_1<1$, then there is $M_0>0$ such that for all $M\ge M_0$ we have the bound
\begin{equation}\label{eq:2.2}
\left\|{\rm Op}_h(a){\rm Op}_h(b)-{\rm Op}_h(E_M(a,b))\right\|_{L^2(Y)\to L^2(Y)}\lesssim h^{M(1-\delta_2-\delta'_1)/2}.
\end{equation}
\end{prop}

Indeed, (\ref{eq:2.2}) is an easy consequence of the inequality (7.17) of \cite{kn:DS} and the bound (\ref{eq:2.1}) (e.g. see Section 4
of \cite{kn:V2}).

Let $y=(y_1,...,y_{d-1})$ be coordinates in $Y$ and let $\eta=(\eta_1,...,\eta_{d-1})$ be the dual variables. 
Let $\mu\in{\bf R}$ be a parameter satisfying $h^{2/3-\epsilon}\le|\mu|\le 1$,
$0<\epsilon\ll 1$ being arbitrary.
Consider in $(0,1)\times Y$ the operator
$$P_0(h,\mu)={\cal D}_t^2+{\cal D}_{y_1}-i\mu+m(t,y,{\cal D}_y;h,\mu)$$
where $t\in (0,1)$, ${\cal D}_t:=-ih\partial_t$, ${\cal D}_{y_1}:=-ih\partial_{y_1}$. We suppose that there are functions
$m_j\in C^\infty([0,1]\times T^*Y)$, $j=0,1,...,$ independent of $h$,  
such that for every integer $\nu\ge 0$ the function $m$ can be written in the form
$$m(t,y,\eta;h,\mu)=\sum_{j=0}^\nu h^jm_j(t,y,\eta;\mu)+h^{\nu+1}m_{\nu+1}^\sharp(t,y,\eta;h,\mu)$$
where $m^\sharp_{\nu+1}\in S^{0,0}_{0,0}(Y)$ uniformly in $t$, $h$ and $\mu$.
We also suppose that ${\rm Im}\,\partial_t^\alpha\partial_y^\beta m_0={\cal O}(|\mu|)$ for all multi-indices $\alpha$ and $\beta$, 
and that for all integers $k,j\ge 0$, $\partial_t^km_j\in S^{0,0}_{0,0}(Y)$ uniformly in $t$ and $\mu$. 
The Taylor expansion at $t=0$ gives, for every integer $\nu\ge 0$,
$$m_j(t,y,\eta;\mu)=\sum_{k=0}^\nu t^km_{k,j}(y,\eta;\mu)+t^{\nu+1}\widetilde m_{\nu+1,j}(t,y,\eta;\mu)$$
where $\widetilde m_{\nu+1,j}\in S^{0,0}_{0,0}(Y)$ uniformly in $t$ and $\mu$. We will write 
$$m\sim\sum_{j=0}^\infty h^jm_j\sim\sum_{k=0}^\infty\sum_{j=0}^\infty  t^kh^jm_{k,j}.$$
Finally, we suppose that the function
$m_{0,0}:=m_0(0,y,\eta;\mu)$ is identically zero on $T^*Y$. In other words, $m={\cal O}(t+h)$ as $t\to 0$. 

Our goal in this section is to build a parametrix, $\widetilde u$, for the solution to the following boundary value
problem:
\begin{equation}\label{eq:2.3}
\left\{
\begin{array}{lll}
 P_0(h,\mu)u(t,y)=0&\mbox{in}& (0,1)\times Y,\\
 u(0,y)=\phi({\cal D}_{y_1})f&\mbox{on}&Y,
\end{array}
\right.
\end{equation}
where $f\in L^2(Y)$ and $\phi\in C_0^\infty({\bf R})$ is independent of $h$ and $\mu$, such that $\phi(\sigma)=1$ for $|\sigma|\le 1$,
$\phi(\sigma)=0$ for $|\sigma|\ge 2$. We will also study the way in which the parametrix depends on the functions $m_{k,j}$. 
We will be looking for $\widetilde u$ in the form
$$\widetilde u(t,y)=(2\pi h)^{-d+1}\int_{T^*Y}e^{-\frac{i}{h}(\langle y-w,\eta\rangle-\varphi(t,y,\eta;\mu))}\Phi_{\epsilon,\delta}(t,\eta_1;\mu)a(t,y,\eta;h,\mu)f(w)dwd\eta$$
where $\Phi_{\epsilon,\delta}=\phi(t/h^\epsilon)\phi(t/|\varrho|^2\delta)$, $0<\delta\ll 1$ being a parameter independent of $h$ and $\mu$ to be fixed later on and the 
function $\varrho$ is defined by
$$\varrho(\eta_1;\mu)=\sqrt{-\eta_1+i\mu},\quad {\rm Im}\,\varrho>0.$$
Clearly, we have $|\varrho|^2\ge|\mu|$. We also have $|\varrho|\le Const$ as long as $\eta_1\in{\rm supp}\,\phi$. 

The phase $\varphi$ is complex-valued, independent of $h$, to be determined later on such that
$\varphi|_{t=0}\equiv 0$. The amplitude $a$ is of the form
$$a=\sum_{j=0}^Mh^ja_j$$
where $M$ is an arbitrary sufficiently large integer and the functions $a_j$ are independent of $h$ to be determined later on such that
$a_0|_{t=0}=\phi(\eta_1)$, $a_j|_{t=0}\equiv 0$ for all $j\ge 1$. Thus we have $\widetilde u(0,y)=\phi({\cal D}_{y_1})f$.
Furthermore, writing formally
$$\widetilde u={\rm Op}_h\left(e^{i\varphi/h}\Phi_{\epsilon,\delta} a\right)f$$
we get
\begin{equation}\label{eq:2.4}
P_0(h,\mu)\widetilde u={\rm Op}_h\left(e^{i\varphi/h}\Phi_{\epsilon,\delta} A_M\right)f+{\rm Op}_h\left(e^{i\varphi/h}A^\sharp_M\right)f+{\cal E}_Mf
\end{equation}
where
$$A_M=\left((\partial_t\varphi)^2+\partial_{y_1}\varphi-\varrho^2-ih\partial_t^2\varphi\right)a-2ih\partial_t\varphi\partial_ta-ih\partial_{y_1}a-h^2
\partial_t^2a$$ $$+e^{-i\varphi/h}E_M(m,e^{i\varphi/h}a),$$
$$A^\sharp_M=\left(-h^2\partial_t^2\Phi_{\epsilon,\delta}-2ih\partial_t\Phi_{\epsilon,\delta}\partial_t\varphi\right)a
-2h^2\partial_t\Phi_{\epsilon,\delta}\partial_ta,$$
$${\cal E}_M={\rm Op}_h\left(m\right){\rm Op}_h\left(e^{i\varphi/h}\Phi_{\epsilon,\delta} a\right)-
{\rm Op}_h\left(E_M(m,e^{i\varphi/h}\Phi_{\epsilon,\delta} a)\right).$$
It is easy to see that $A^\sharp_M=0$ for $0\le t\le\min\{h^\epsilon,\delta|\rho|^2\}$. We will now expand the function $A_M$ in powers of $h$. The most difficult
is to expand the last term. To do so 
we need the following

\begin{lemma} For every multi-index $\beta$ we have the identity
\begin{equation}\label{eq:2.5}
\frac{(-i)^{|\beta|}}{|\beta|!}\partial_y^\beta\left(e^{i\varphi/h}\right)=e^{i\varphi/h}\sum_{k=0}^{|\beta|}h^{-k}G_{k}^{(\beta)}(\varphi)
\end{equation}
where the functions $G_{k}^{(\beta)}$ are independent of $h$, $G_{0}^{(0)}=1$, $G_{0}^{(\beta)}=0$ for $|\beta|\ge 1$. For $k\ge 1$, $|\beta|\ge 1$ they are of
the form
$$G_{k}^{(\beta)}(\varphi)=\sum_{|\gamma_j|\ge 1,\,|\gamma_1|+...+|\gamma_k|=|\beta|}c_{\gamma_1,...,\gamma_k;k,\beta}
\prod_{j=1}^k\partial_y^{\gamma_j}\varphi$$
where the coefficients $c_{\gamma_1,...,\gamma_k;k,\beta}$ are constants. In particular, if $\beta=(\beta_1,...,\beta_{d-1})$,
then
$$G_{|\beta|}^{(\beta)}(\varphi)=\frac{1}{|\beta|!}\prod_{j=1}^{d-1}
\left(\frac{\partial\varphi}{\partial y_j}\right)^{\beta_j}.$$
\end{lemma}

This lemma can be easilly proved by induction in $|\beta|$ and therefore we omit the details. Using (\ref{eq:2.5}) we can write
$$\frac{(-i)^{|\alpha|}}{|\alpha|!}\partial_y^\alpha\left(e^{i\varphi/h}a\right)=\sum_{0\le|\beta|\le|\alpha|}
c_{\alpha,\beta}\frac{(-i)^{|\alpha-\beta|}}{|\alpha-\beta|!}\partial_y^{\alpha-\beta}\left(e^{i\varphi/h}\right)\partial_y^{\beta}a$$
$$=e^{i\varphi/h}\sum_{0\le|\beta|\le|\alpha|}\sum_{k=0}^{|\alpha|-|\beta|}c_{\alpha,\beta}h^{-k}G_{k}^{(\alpha-\beta)}(\varphi)
\partial_y^{\beta} a$$
$$=e^{i\varphi/h}h^{-|\alpha|}\sum_{\nu=0}^{|\alpha|}h^\nu\sum_{|\beta|\le\nu}c_{\alpha,\beta}
G_{|\alpha|-\nu}^{(\alpha-\beta)}(\varphi)\partial_y^{\beta} a$$
where the coefficients $c_{\alpha,\beta}$ are constants. Thus we get the expansion
$$e^{-i\varphi/h}\frac{(-ih)^{|\alpha|}}{|\alpha|!}\partial_y^\alpha\left(e^{i\varphi/h}a\right)=G_{|\alpha|}^{(\alpha)}(\varphi)a+\sum_{\nu=1}^{M+|\alpha|} h^\nu a_\nu^{(\alpha)}$$
where 
$$a_\nu^{(\alpha)}=\sum_{\nu'=(\nu-|\alpha|)_+}^{\nu-1}\sum_{|\beta|\le\nu-\nu'}
c_{\alpha,\beta}G_{|\alpha|-\nu+\nu'}^{(\alpha-\beta)}(\varphi)
\partial_y^{\beta} a_{\nu'}$$
where $b_+:=\max{\{0,b}\}$.
 Clearly, we have
$a_\nu^{(0)}=0$. We now expand the function $m$ as
$$m=\sum_{j=0}^{M+1} h^jm_j+h^{M+2}m^\sharp_{M+2}.$$
Using the above identities we can write
$$e^{-i\varphi/h}E_M(m,e^{i\varphi/h}a)=g_M(m_0,\varphi)a+\sum_{j=0}^{M} h^{j+1} E_j^{(M)}+h^{M+2}\widetilde E_{M+2}$$
where 
$$g_M(m_0,\varphi)=\sum_{0\le|\alpha|\le M}\partial_\eta^\alpha m_0G_{|\alpha|}^{(\alpha)}(\varphi),$$
$$E_j^{(M)}=\sum_{0\le|\alpha|\le M}\sum_{\ell=0}^{j} \partial_\eta^\alpha m_\ell a^{(\alpha)}_{j+1-\ell}
+\sum_{0\le|\alpha|\le M}G_{|\alpha|}^{(\alpha)}(\varphi)\sum_{\ell=1}^{j+1} \partial_\eta^\alpha m_\ell a_{j+1-\ell},$$
$$\widetilde E_{M+2}=\sum_{0\le|\alpha|\le M}\partial_\eta^\alpha m_{M+2}^\sharp\left(
G_{|\alpha|}^{(\alpha)}(\varphi)a+\sum_{\nu=1}^{M+|\alpha|} h^\nu a_\nu^{(\alpha)}\right)$$
$$+\sum_{0\le|\alpha|\le M}\sum_{\nu=M+1}^{2M+|\alpha|}h^{\nu-M-1}\sum_{\ell=0}^{\nu+1} \partial_\eta^\alpha m_\ell a^{(\alpha)}_{\nu+2-\ell}.$$
We can expand the function $A_M$ as follows:
$$A_M=\left((\partial_t\varphi)^2+\partial_{y_1}\varphi-\varrho^2+g_M(m_0,\varphi)\right)a$$
$$-\sum_{j=0}^Mh^{j+1}\left(2i\partial_t\varphi\partial_t a_j+i\partial_{y_1}a_j+i\partial_t^2\varphi a_j+\partial_t^2a_{j-1}-E_{j}^{(M)}\right)$$
$$-h^{M+2}\partial_t^2a_M+h^{M+2}\widetilde E_{M+2}.$$
We would like to determine the functions $\varphi$ and $a_j$, $j=0,1,...,M$, so that $A_M={\cal O}(t^{M})+{\cal O}(h^{M+2})$. To this end, we let
the function $\varphi$ satisfy the following eikonal equation mod ${\cal O}(t^{M})$:
\begin{equation}\label{eq:2.6}
(\partial_t\varphi)^2+\partial_{y_1}\varphi-\rho^2+g_M(m_0,\varphi)=t^{M}R_M
\end{equation}
where the function $R_M$ is smooth up to $t=0$. The functions $a_j$ satisfy the transport equations mod ${\cal O}(t^{M})$:
\begin{equation}\label{eq:2.7}
2i\partial_t\varphi\partial_t a_j+i\partial_{y_1}a_j+i\partial_t^2\varphi a_j+\partial_t^2a_{j-1}-E_{j}^{(M)}=t^{M}Q_M^{(j)},\quad
0\le j\le M,
\end{equation}
where $a_{-1}=0$ and the functions $Q_M^{(j)}$ are smooth up to $t=0$. Then we get
\begin{equation}\label{eq:2.8}
A_M=t^MB_M+h^MC_M
\end{equation}
where
$$B_M=R_Ma-\sum_{j=0}^Mh^{j+1}Q_M^{(j)},$$
$$C_M=-h^{2}\partial_t^2a_M+h^{2}\widetilde E_{M+2}.$$
We will first solve equation (\ref{eq:2.6}). We will be looking for $\varphi$
in the form
$$\varphi=\sum_{k=1}^M t^k\varphi_k$$
with functions $\varphi_k$ independent of $t$. We have
$$(\partial_t\varphi)^2=\sum_{K=0}^{2M-2}t^K\sum_{k+j=K}(k+1)(j+1)\varphi_{k+1}\varphi_{j+1},$$
$$\partial_{y_1}\varphi=\sum_{K=1}^M t^K\partial_{y_1}\varphi_K,$$
$$G_{|\alpha|}^{(\alpha)}(\varphi)=\sum_{k=|\alpha|}^{M|\alpha|}t^k{\cal V}_k^{(\alpha)}(\varphi_1,...,\varphi_k),$$
where ${\cal V}_0^{(0)}=1$, while for $|\alpha|\ge 1$, $k\ge|\alpha|$ we have 
$${\cal V}_k^{(\alpha)}=\sum_{k_j\ge 1,\,k_1+...+k_{|\alpha|}=k}\sum_{|\gamma_1|=...=|\gamma_{|\alpha|}|=1}c_{k_1,...,k_{|\alpha|},\gamma_1,...,\gamma_{|\alpha|}}\prod_{j=1}^{|\alpha|}
\partial_y^{\gamma_j}\varphi_{k_j}$$
where the coefficients $c_{k_1,...,k_{|\alpha|},\gamma_1,...,\gamma_{|\alpha|}}$ are real constants. We now expand the function $m_0$ as
$$m_0=\sum_{k=1}^{M-1} t^k m_{k,0}+t^{M}\widetilde m_{M,0}$$
where the functions $m_{k,0}$, $k=1,...,M-1,$ are independent of $t$, while the function $\widetilde m_{M,0}$ is smooth up to $t=0$. Hence
$$g_M(m_0,\varphi)=\sum_{K=1}^{M-1} t^K \sum_{0\le|\alpha|\le K-1}\sum_{\nu=1}^{K-|\alpha|}\partial_\eta^\alpha m_{\nu,0}{\cal V}_{K-\nu}^{(\alpha)}+t^{M}\widetilde g_M$$
where
$$\widetilde g_M=\sum_{0\le|\alpha|\le M}\partial_\eta^\alpha \widetilde m_{M,0}G_{|\alpha|}^{(\alpha)}(\varphi)
+\sum_{0\le|\alpha|\le M}t^{-1}\partial_\eta^\alpha m_0\sum_{k=M-1}^{M|\alpha|}t^{k-M+1}{\cal V}_k^{(\alpha)}$$
$$+\sum_{K=M}^{2M-1} t^{K-M} \sum_{0\le|\alpha|\le K-1}\sum_{\nu=1}^{K-|\alpha|}\partial_\eta^\alpha m_{\nu,0}{\cal V}_{K-\nu}^{(\alpha)}.$$
Inserting the above identities into equation (\ref{eq:2.6}) and comparing the coefficients of all powers $t^K$, $0\le K\le M-1$, we get the following relations for the functions $\varphi_k$:
$$\varphi_1^2=\varrho^2,$$
\begin{equation}\label{eq:2.9}
\sum_{k+j=K}(k+1)(j+1)\varphi_{k+1}\varphi_{j+1}+\partial_{y_1}\varphi_K+m_{K,0}$$ $$
+\sum_{1\le|\alpha|\le K-1}\sum_{\nu=1}^{K-|\alpha|}\partial_\eta^\alpha m_{\nu,0}{\cal V}_{K-\nu}^{(\alpha)}(\varphi_1,...,\varphi_{K-\nu})=0,
\end{equation}
$1\le K\le M-1$, where the second sum is zero if $K=1$. Then equation (\ref{eq:2.6}) is satisfied with
$$R_M=\widetilde g_M+\sum_{K=M}^{2M-2}t^{K-M}\sum_{k+j=K}(k+1)(j+1)\varphi_{k+1}\varphi_{j+1}+\partial_{y_1}\varphi_M.$$
We take $\varphi_1=\varrho$. Then, given $\varphi_1,...,\varphi_K$, we can find $\varphi_{K+1}$ from equation (\ref{eq:2.9}) in a unique way.
Thus we can find all functions $\varphi_1,...,\varphi_M$. In what follows in this section, given a function $p>0$ on $T^*Y$ and $k\in{\bf R}$, we will denote by $S^k(p)$ the set of all functions $a\in C^\infty(T^*Y)$ satisfying $\partial_y^\alpha a={\cal O}_\alpha(p^k)$ for all multi-indices 
$\alpha$. 

\begin{lemma} For all integers $k\ge 1$ we have $\varphi_k\in S^{3-2k}(|\varrho|)$ and 
\begin{equation}\label{eq:2.10}
\left|{\rm Im}\,\partial_y^\alpha\varphi_k\right|\le C_{k,\alpha}|\varrho|^{2-2k}{\rm Im}\,\varrho.
\end{equation}
Moreover, if $0<t\le 2\delta|\varrho|^2$ with a constant $\delta>0$ small enough, we have
\begin{equation}\label{eq:2.11}
{\rm Im}\,\varphi\ge t{\rm Im}\,\varrho/2.
\end{equation}
We also have that the functions $\varphi_k$ and $\varphi_{k+1}+\frac{m_{k,0}}{2(k+1)\varrho}$ are independent of all $m_{\ell,0}$ with $\ell\ge k$.
\end{lemma}

{\it Proof.} We will proceed by induction in $k$. Clearly, $\varphi_1\in S^1(|\varrho|)$. Suppose now that
$\varphi_k\in S^{3-2k}(|\varrho|)$ for $1\le k\le K$. This implies 
$$\varphi_{k+1}\varphi_{K-k+1}\in S^{2-2K}(|\varrho|),\quad 1\le k\le K-1,$$
$$\partial_{y_1}\varphi_K\in S^{3-2K}(|\varrho|),$$
$${\cal V}_{K-\nu}^{(\alpha)}(\varphi_1,...,\varphi_{K-\nu})\in S^{3|\alpha|+2\nu-2K}(|\varrho|).$$
Therefore, by equation (\ref{eq:2.9}) we conclude that $\varrho\varphi_{K+1}\in S^{2-2K}(|\varrho|)$, which implies
$\varphi_{K+1}\in S^{1-2K}(|\varrho|)$ as desired. The bound (\ref{eq:2.10}) can also been easily proved by induction in $k$. Indeed,
differentiating equation (\ref{eq:2.9}) allows to express $\partial_y^\alpha\varphi_{K+1}$ in terms of the functions $m_{1,0},...,
m_{K,0}$, $\varphi_1,...,\varphi_K$ and their derivatives. Thus we can bound $|{\rm Im}\,\partial_y^\alpha\varphi_{K+1}|$ by using 
that  by assumption we have $${\rm Im}\,\partial_y^\alpha m_{k,0}={\cal O}(|\mu|)={\cal O}(|\varrho|{\rm Im}\,\varrho)$$  together with the inequality
$$|{\rm Im}\,(z_1...z_k)|\le C_k|z_1|...|z_k|\sum_{j=1}^k\frac{|{\rm Im}\,z_j|}{|z_j|}.$$
The bound (\ref{eq:2.11}) follows from (\ref{eq:2.10}). We have, for $0<t\le 2\delta|\varrho|^2$,
$${\rm Im}\,\varphi=\sum_{k=1}^Mt^k{\rm Im}\,\varphi_k\ge t{\rm Im}\,\varrho\left(1-C\sum_{k=1}^{M-1}t^k|\varrho|^{-2k}\right)$$
$$\ge t{\rm Im}\,\varrho(1-{\cal O}(\delta))\ge t{\rm Im}\,\varrho/2$$
provided $\delta$ is taken small enough. The last assertion also follows by induction in $k$. Indeed, by equation (\ref{eq:2.9}) we
can express the function $\varphi_{K+1}+\frac{m_{K,0}}{2(K+1)\varrho}$ in terms of the functions $m_{1,0},...,
m_{K-1,0}$, $\varphi_1,...,\varphi_K$ and their derivatives. Therefore, it is independent of all $m_{\ell,0}$ with $\ell\ge K$,
provided so are $\varphi_k$ with $k\le K$. 
\eproof

The above lemma implies the following

\begin{lemma} For all $|\beta|\ge 0$, $0\le k\le |\beta|$, we have the identity
\begin{equation}\label{eq:2.12}
G_{k}^{(\beta)}(\varphi)=\sum_{\nu=k}^{kM}t^\nu\Theta_\nu^{(k,\beta)}(m_{1,0},...,m_{\nu-k,0})
\end{equation}
where the function $\Theta_\nu^{(k,\beta)}$ is independent of $t$ and all $m_{\ell,0}$ with $\ell\ge \nu-k+1$. Moreover,
$\Theta_0^{(0,0)}=1$, $\Theta_\nu^{(0,\beta)}=0$ for all $\nu\ge 0$ if $|\beta|\ge 1$, while for 
all $|\beta|\ge 1$, $1\le k\le |\beta|$,
we have $\Theta_\nu^{(k,\beta)}\in S^{-2\nu}(|\varrho|)$. 
\end{lemma}

{\it Proof.} It is clear from the definition of the function $G_{k}^{(\beta)}(\varphi)$ that we have (\ref{eq:2.12}) with $\Theta_\nu^{(k,\beta)}$ being a linear 
combination of functions of the form
$$\prod_{j=1}^k\partial_y^{\gamma_j}\varphi_{\nu_j}$$
where $1\le\nu_j\le M$ are integers such that $\nu_1+...+\nu_k=\nu$. Hence $\nu$ must satisfy $k\le\nu\le kM$
and each $\nu_j$
must satisfy $\nu_j\le \nu-k+1$. Therefore, the lemma is an immediate consequence of Lemma 2.3.
\eproof

We will now be looking for solutions to the equations (\ref{eq:2.7}) in the form
$$a_j=\sum_{k=0}^Mt^ka_{k,j}.$$
We have
$$\partial_t\varphi\partial_t a_j=\sum_{k=0}^{2M-2}t^k\sum_{\nu=0}^{k}(k-\nu+1)(\nu+1)\varphi_{k-\nu+1} a_{\nu+1,j},$$
$$\partial_{y_1}a_j=\sum_{k=0}^Mt^k\partial_{y_1}a_{k,j},$$
$$\partial_t^2\varphi a_j=\sum_{k=0}^{2M-2}t^k\sum_{\nu=0}^{k}(k-\nu+1)(k-\nu+2)\varphi_{k-\nu+2} a_{\nu,j},$$
$$\partial_t^2a_{j-1}=\sum_{k=0}^{M-2}t^k(k+1)(k+2)a_{k+2,j-1}.$$
We will now use Lemma 2.4 to prove the following

\begin{lemma} We have the identity
$$E_{j}^{(M)}=\sum_{k=0}^{M(M+2)}t^kE_{k,j}^{(M)}$$
with functions $E_{k,j}^{(M)}$ independent of $t$ having the form
$$E_{k,j}^{(M)}=\sum_{\nu=0}^k\sum_{\ell=0}^j\sum_{|\gamma|\le M}\Psi_{\nu,\ell,\gamma}^{(k,j)}\partial_y^\gamma a_{\nu,\ell}$$
where the function $\Psi_{\nu,\ell,\gamma}^{(k,j)}\in S^{-2(k-\nu)}(|\varrho|)$ depends only on $m_{k',j'}$
with $k'\le k$, $j'\le j+1$.
\end{lemma}

{\it Proof.} Using (\ref{eq:2.12}) we expand the function $a_\nu^{(\alpha)}$ as
$$a_\nu^{(\alpha)}=\sum_{k=0}^{M(|\alpha|+1)}t^ka_{k,\nu}^{(\alpha)}$$
where
$$a_{k,\nu}^{(\alpha)}=\sum_{\nu'=(\nu-|\alpha|)_+}^{\nu-1}\sum_{k'=0}^{k+\nu-\nu'-|\alpha|}\sum_{|\beta|\le\nu-\nu'}
 c_{\alpha,\beta}\Theta_{k-k'}^{(|\alpha|-\nu+\nu',\alpha-\beta)}
\partial_y^{\beta} a_{k',\nu'}.$$
Hence
$$\sum_{\ell=0}^{j} \partial_\eta^\alpha m_\ell a^{(\alpha)}_{j+1-\ell}=\sum_{k=0}^{M(|\alpha|+2)}t^k\sum_{k_1=0}^k
\sum_{j_1=1}^{j+1}\partial_\eta^\alpha m_{k-k_1,j+1-j_1}a^{(\alpha)}_{k_1,j_1}$$
$$=\sum_{k=0}^{M(|\alpha|+2)}t^k\sum_{k'=0}^k\sum_{j'=0}^j\sum_{|\beta|\le |\alpha|}F_{k',j',\beta}^{(k,j,\alpha)}\partial_y^{\beta} a_{k',j'}$$
where
$$F_{k',j',\beta}^{(k,j,\alpha)}=c_{\alpha,\beta}\sum_{(k_1,j_1)\in \Xi(k,j,k',j',\alpha,\beta)}\partial_\eta^\alpha m_{k-k_1,j+1-j_1} 
\Theta_{k_1-k'}^{(|\alpha|-j_1+j',\alpha-\beta)}$$
where $\Xi(k,j,k',j',\alpha,\beta)$ denotes the set of all integers $(k_1,j_1)$ satisfying $0\le k_1\le k$, $1\le j_1\le j+1$,
$j'+|\beta|\le j_1\le j'+|\alpha|$, $k_1+j_1\ge k'+j'+|\alpha|$. Clearly, if this set is empty, the above sum is zero.
By Lemma 2.4 we get
$$F_{k',j',\beta}^{(k,j,\alpha)}\in S^{-2(k-k')}(|\varrho|).$$
Furthermore, we have
$$\sum_{\ell=1}^{j+1} \partial_\eta^\alpha m_\ell a_{j+1-\ell}=\sum_{k=0}^{2M}t^k\sum_{k'=0}^k\sum_{j'=0}^{j}
\partial_\eta^\alpha m_{k-k',j+1-j'}a_{k',j'}.$$
Hence
$$G_{|\alpha|}^{(\alpha)}(\varphi)\sum_{\ell=1}^{j+1} \partial_\eta^\alpha m_\ell a_{j+1-\ell}=
\sum_{k=0}^{M(|\alpha|+2)}t^k\sum_{k_1=0}^{k-|\alpha|}\sum_{k'=0}^{k_1}\sum_{j'=0}^{j}{\cal V}_{k-k_1}^{(\alpha)}
\partial_\eta^\alpha m_{k_1-k',j+1-j'}a_{k',j'}$$
$$=\sum_{k=0}^{M(|\alpha|+2)}t^k\sum_{k'=0}^{k-|\alpha|}\sum_{j'=0}^{j}\widetilde F_{k',j'}^{(k,j,\alpha)}a_{k',j'}$$
where
$$\widetilde F_{k',j'}^{(k,j,\alpha)}=\sum_{k_1=k'}^{k-|\alpha|}{\cal V}_{k-k_1}^{(\alpha)}
\partial_\eta^\alpha m_{k_1-k',j+1-j'}\in S^{-2(k-k')}(|\varrho|).$$
It follows from the above identities that the desired expansion of the function $E_{j}^{(M)}$ holds with
$$\Psi_{\nu,\ell,\gamma}^{(k,j)}=\sum_{|\gamma|\le|\alpha|\le M}F_{\nu,\ell,\gamma}^{(k,j,\alpha)}\in S^{-2(k-\nu)}(|\varrho|),\quad|\gamma|\ge 1,$$
$$\Psi_{\nu,\ell,0}^{(k,j)}=\sum_{0\le|\alpha|\le M}\left(F_{\nu,\ell,0}^{(k,j,\alpha)}+\widetilde F_{\nu,\ell}^{(k,j,\alpha)}\right)\in S^{-2(k-\nu)}(|\varrho|).$$
\eproof

Inserting the above identities into equations (\ref{eq:2.7}) and comparing the coefficients of all powers $t^k$, $0\le k\le M-1$, we get the following
relations for the functions $a_{k,j}$:
\begin{equation}\label{eq:2.13}
2i\sum_{\nu=0}^{k}(k-\nu+1)(\nu+1)\varphi_{k-\nu+1} a_{\nu+1,j}+
i\sum_{\nu=0}^{k}(k-\nu+1)(k-\nu+2)\varphi_{k-\nu+2} a_{\nu,j}$$ $$+i\partial_{y_1}a_{k,j}+(k+1)(k+2)a_{k+2,j-1}=E_{k,j}^{(M)}
\end{equation}
 and $a_{0,0}=\phi(\eta_1)$, $a_{0,j}=0$, $j\ge 1$, $a_{k,-1}=0$, $k\ge 0$. Since the function $E_{k,j}^{(M)}$ depends only on $a_{k',j'}$ with $k'\le k$, $j'\le j$, it is easy to see that
we can determine all $a_{k,j}$ from equations (\ref{eq:2.13}). Observe also that ${\rm supp}_{\eta_1}a_{k,j}\equiv{\rm supp}\,\phi(\eta_1)$. 
Then the equations (\ref{eq:2.7}) are satisfied with 
$$Q_M^{(j)}=2i\sum_{k=M}^{2M-2}t^{k-M}\sum_{\nu=0}^{M-1}(k-\nu+1)(\nu+1)\varphi_{k-\nu+1} a_{\nu+1,j}+
i\partial_{y_1}a_{M,j}$$
$$+i\sum_{k=M}^{2M-2}t^{k-M}\sum_{\nu=0}^{M}(k-\nu+1)(k-\nu+2)\varphi_{k-\nu+2} a_{\nu,j}$$ 
$$-\sum_{k=M}^{M(M+2)}t^{k-M}\sum_{\nu=0}^M\sum_{\ell=0}^j\sum_{|\gamma|\le M}\Psi_{\nu,\ell,\gamma}^{(k,j)}\partial_y^\gamma a_{\nu,\ell}.$$
We will now prove the following

\begin{lemma} For all $k\ge 0,\,j\ge 0$, we have $a_{k,j}\in S^{-2k-3j}(|\varrho|)$. Moreover, the function
$$a_{k,j}+\frac{(k+j)!}{k!}\frac{\phi m_{k+j,0}}{(-2i\varrho)^{j+2}}$$
depends only on $m_{\nu,0}$ with $\nu\le k+j-1$ and $m_{\nu,\ell}$ with $1\le\ell\le j+1$, $\nu+\ell\le k+j$.
\end{lemma}

{\it Proof.} We will proceed by induction. Clearly, the first assertion is trivial for $k=0$ and all $j\ge 0$. Suppose that it is true
for all $j\le J-1$, $k\ge 0$, and for $j=J$, $k\le K$. We have to show that it is true for $j=J$ and $k=K+1$. In view of Lemmas 2.3 and 2.5
we have
$$\varphi_{k-\nu+1} a_{\nu+1,j}\in S^{-2k-3j-1}(|\varrho|),\quad 0\le\nu\le k-1,k=K,j=J,$$
$$\varphi_{k-\nu+2} a_{\nu,j}\in S^{-2k-3j-1}(|\varrho|),\quad 0\le\nu\le k,k=K,j=J,$$
$$\partial_{y_1}a_{k,j}\in S^{-2k-3j}(|\varrho|),\quad k=K,j=J,$$
$$a_{k+2,j-1}\in S^{-2k-3j-1}(|\varrho|),\quad k=K,j=J,$$
$$E_{k,j}^{(M)}\in S^{-2k-3j}(|\varrho|),\quad k=K,j=J.$$
Therefore, by equation (\ref{eq:2.13}) we conclude $\varrho a_{K+1,J}\in S^{-2K-3J-1}(|\varrho|)$, which implies
$a_{K+1,J}\in S^{-2K-3J-2}(|\varrho|)$ as desired. 

We turn now to the proof of the second assertion. We will first prove it for $j=0$ and all $k\ge 0$ by induction in $k$. 
It is trivially fulfilled for $k=0$. Suppose it is fulfilled for all $k\le K$. We have to show it is fulfilled for $k=K+1$.
By Lemma 2.5 we have that $E_{K,0}^{(M)}$ depends on 
$m_{\nu,\ell}$ with $0\le\ell\le 1$, $\nu\le K$. Therefore, by equation (\ref{eq:2.13}) with $j=0$, $k=K$, we get that
$$2i(K+1)\varrho a_{K+1,0}+i(K+1)(K+2)\varphi_{K+2}a_{0,0}$$
depends on $m_{\nu,\ell}$ with $0\le\ell\le 1$, $0\le\nu\le K$. We now use the last assertion in Lemma 2.3 to conclude that
$$a_{K+1,0}-(2\varrho)^{-2}m_{K+1,0}\phi$$
depends on $m_{\nu,\ell}$ with $0\le\ell\le 1$, $0\le\nu\le K$, as desired.

We will now proceed by induction in both $j$ and $k$. The assertion is trivially fulfilled for $k=0$. Suppose it is true
for all $j\le J-1$, $k\ge 0$ and for $j=J$, $k\le K$ with some integers $K\ge 0$, $J\ge 1$. We will prove it for $j=J$ and $k=K+1$. 
To this end we will use 
equation (\ref{eq:2.13}) with $k=K$ and $j=J$. By Lemma 2.5 we have that $E_{K,J}^{(M)}$ depends on 
$m_{\nu,0}$ with $\nu\le K+J$ and $m_{\nu,\ell}$ with $1\le\ell\le J+1$, $\nu+\ell\le K+J$. Therefore, we get that so does the function
$$2i\varrho(K+1)a_{K+1,J}+(K+1)(K+2)a_{K+2,J-1}.$$
On the other hand, since the assertion is supposed to be fulfilled for $j=J-1$ and all $k$, we have that the function
$$(K+2)a_{K+2,J-1}+\frac{(K+1+J)!}{(K+1)!}\frac{\phi m_{K+1+J,0}}{(-2i\varrho)^{J+1}}$$
depends only on $m_{\nu,0}$ with $\nu\le K+J$ and $m_{\nu,\ell}$ with $1\le\ell\le J$, $\nu+\ell\le K+1+J$.
Thus we conclude that the function
$$a_{K+1,J}+\frac{(K+1+J)!}{(K+1)!}\frac{\phi m_{K+1+J,0}}{(-2i\varrho)^{J+2}}$$
depends only on $m_{\nu,0}$ with $\nu\le K+J$ and $m_{\nu,\ell}$ with $1\le\ell\le J+1$, $\nu+\ell\le K+1+J$, as desired.
\eproof

It is clear from the equations (\ref{eq:2.13}) that the functions $a_{k,j}$ are well-defined for all $\mu\neq 0$ because so is
the function $\varrho^{-1}$. The condition $|\mu|\ge h^{2/3-\epsilon}$ is only used to show that the above construction provides
a parametrix for the boundary value problem (\ref{eq:2.3}) (see the proof of Proposition 2.8 below). It is also possible to bound the derivatives
of $a_{k,j}$ with repspect to the variable $\eta$ uniformly in $\mu$. Thus, although we do not need this information in the analysis
that follows, for some values of $\mu$ we can describe completely the class
of symbols the functions $a_{k,j}$ belong to. For example, we have the following

\begin{lemma} Let $|\mu|\ge h^{2/3}$. Then we have $a_{k,j}\in S_{0,2/3}^{2k/3+j,0}$.
\end{lemma}

{\it Proof.} In the same way as above one can show that the functions $a_{k,j}$ satisfy the bounds
$$\left|\partial_y^\alpha\partial_\eta^\beta a_{k,j}\right|\le C_{\alpha,\beta}|\varrho|^{-2k-3j-2\beta_1}
\le C_{\alpha,\beta}|\mu|^{-k-3j/2-\beta_1}$$
for all multi-indices $\alpha$ and $\beta=(\beta_1,...)$, which clearly imply the lemma.
\eproof

We will now show that $\widetilde u$ provides the desired parametrix. Recall that $\widetilde u$ depends on the
parameter $M$. We have the following

\begin{prop} Let $h^{2/3-\epsilon}\le|\mu|\le 1$. Then there is $M_0>0$ depending on $\epsilon$ such that for
$M\ge M_0$ we have the bounds
\begin{equation}\label{eq:2.14}
\left\|\widetilde u\right\|_{L^2((0,1)\times Y)}\lesssim h^{-d}\|f\|_{L^2(Y)},
\end{equation}
\begin{equation}\label{eq:2.15}
\left\|P_0(h,\mu)\widetilde u\right\|_{L^2((0,1)\times Y)}\lesssim h^{\epsilon M/2}\|f\|_{L^2(Y)}.
\end{equation}
\end{prop}

{\it Proof.} We will use the identity (2.4). Observe first that by Lemma 2.3 we have $\partial_t\varphi,\varphi/t\in S^1(|\varrho|)$ and
$${\rm Im}\,\varphi\ge t{\rm Im}\,\varrho/2\ge \frac{t|\mu|}{4|\varrho|}$$
as long as $0<t\le 2\delta|\varrho|^2$. Thus, by Lemma 2.2 we get
\begin{equation}\label{eq:2.16}
\left|\partial_y^\alpha\left(e^{i\varphi/h}\right)\right|\lesssim \left(\frac{t}{h}\right)^{|\alpha|}\exp\left(-\frac{t|\mu|}{4h|\varrho|}\right).
\end{equation}
In particular (\ref{eq:2.16}) implies
\begin{equation}\label{eq:2.17}
\left|\partial_y^\alpha\left(e^{i\varphi/h}\right)\right|\lesssim \exp\left(-h^{-3\epsilon/2}/5\right)
\end{equation}
for $h^\epsilon\le t\le 2h^\epsilon$ or $\delta|\varrho|^2\le t\le 2\delta|\varrho|^2$.
Next by Lemma 2.6 we have, with any $\nu\ge 0$,
$$t^{k-\nu}a_{k,j}\in S^{-2\nu-3j}(|\varrho|)\subset S^{-\nu-3j/2}(|\mu|)\subset S^{-2\nu/3-j}(h)$$
as long as $0<t\le 2\delta|\varrho|^2$. Hence
\begin{equation}\label{eq:2.18}
t^{k-\nu}h^ja_{k,j}\in S^{-2\nu/3}(h)
\end{equation}
which implies
\begin{equation}\label{eq:2.19}
\partial_t^\nu a\in S^{-2\nu/3}(h),\quad\nu\ge 0.
\end{equation}
It follows from (\ref{eq:2.19}) that $\partial_y^\alpha A_M^\sharp={\cal O}_\alpha(1)$ which together with 
(\ref{eq:2.17}) imply
$$\left|\partial_y^\alpha\left(e^{i\varphi/h} A_M^\sharp\right)\right|\lesssim \exp\left(-h^{-3\epsilon/2}/5\right).$$
This bound together with (\ref{eq:2.1}) yield
\begin{equation}\label{eq:2.20}
\left\|{\rm Op}_h\left(e^{i\varphi/h} A_M^\sharp\right)\right\|_{L^2(Y)\to L^2(Y)}\lesssim \exp\left(-h^{-3\epsilon/2}/5\right).
\end{equation}
By (\ref{eq:2.16}) and (\ref{eq:2.19}) with $\nu=0$ we also have
$$\left|\partial_y^\alpha\left(e^{i\varphi/h}\Phi_{\epsilon,\delta}a\right)\right|\lesssim  h^{-(1-\epsilon)|\alpha|}.$$
Therefore by Proposition 2.1 we get
\begin{equation}\label{eq:2.21}
\left\|{\cal E}_M\right\|_{L^2(Y)\to L^2(Y)}\lesssim h^{\epsilon M/2}
\end{equation}
for $M$ big enough, while the bound (\ref{eq:2.1}) yields
\begin{equation}\label{eq:2.22}
\left\|{\rm Op}_h\left(e^{i\varphi/h}\Phi_{\epsilon,\delta}a \right)\right\|_{L^2(Y)\to L^2(Y)}\lesssim h^{-d}.
\end{equation}
 To bound the norm of the first operator in the right-hand side of (\ref{eq:2.4}) we will make use of 
(\ref{eq:2.8}). Using (\ref{eq:2.18}) together with Lemmas 2.3 and 2.5 it is not hard to check that the functions $B_M$ and
$C_M$ belong to the spaces $S^{-2M}(|\varrho|)$ and $S^{-3M}(|\varrho|)$, respectively, uniformly in $h$ and in 
$t\in{\rm supp}\,\Phi_{\epsilon,\delta}$. This fact together with (\ref{eq:2.16}) lead to the bounds
$$\left|\partial_y^\alpha\left(e^{i\varphi/h}\Phi_{\epsilon,\delta}A_M\right)\right|\le
t^M\left|\partial_y^\alpha\left(e^{i\varphi/h}\Phi_{\epsilon,\delta}B_M\right)\right|+h^M\left|\partial_y^\alpha\left(e^{i\varphi/h}\Phi_{\epsilon,\delta}C_M\right)\right|$$
$$\lesssim h^{-|\alpha|}\left(t^M|\varrho|^{-2M}e^{-\frac{t|\mu|}{4h|\varrho|}}+h^M|\varrho|^{-3M}\right)$$ $$
\lesssim  h^{M-|\alpha|}\left(|\varrho|^{-M}|\mu|^{-M}+|\varrho|^{-3M}\right)$$ $$\lesssim  h^{M-|\alpha|}|\mu|^{-3M/2}
\lesssim  h^{3\epsilon M/2-|\alpha|}.$$
Thus by (\ref{eq:2.1}) we get
\begin{equation}\label{eq:2.23}
\left\|{\rm Op}_h\left(e^{i\varphi/h}\Phi_{\epsilon,\delta} A_M\right)\right\|_{L^2(Y)\to L^2(Y)}\lesssim h^{\epsilon M}.
\end{equation}
Combining (\ref{eq:2.20}), (\ref{eq:2.21}) and (\ref{eq:2.23}) we obtain, for $M$ big enough,
\begin{equation}\label{eq:2.24}
\left\|P_0(h,\mu)\widetilde u\right\|_{L^2(Y)}^2\lesssim h^{\epsilon M}\|f\|_{L^2(Y)}^2
\end{equation}
uniformly in $t$. Integrating (\ref{eq:2.24}) with respect to $t$ we get (\ref{eq:2.15}). The bound (\ref{eq:2.14})
follows in the same way from (\ref{eq:2.22}).
\eproof

Define the operator $\widetilde N(h,\mu)$ by $\widetilde Nf:={\cal D}_t\widetilde u|_{t=0}$. We have
$$
\widetilde N={\rm Op}_h\left(a\partial_t\varphi|_{t=0}-ih\partial_ta|_{t=0}\right)={\rm Op}_h\left(\varrho\phi-i\sum_{j=0}^Mh^{j+1}a_{1,j}\right).
$$
Let $s,k\ge 0$ be arbitrary integers such that $k\le 3s+2$ and take $M\gg s$. Set
$$\widetilde N_{s,k}={\rm Op}_h\left(\varrho^{k+1}\phi-i\sum_{j=0}^{s-1}h^{j+1}\varrho^ka_{1,j}\right)$$
where the sum is zero if $s=0$. Let $\phi_1\in C_0^\infty$ be such that $\phi_1=1$ on supp$\,\phi$. Clearly, $(1-\phi_1(\eta_1))a_{1,j}\equiv 0$. 

\begin{lemma} For $|\mu|\ge h^{2/3}$ we have the estimate
\begin{equation}\label{eq:2.25}
\left\|\widetilde N(h,\mu){\rm Op}_h\left(\phi_1\varrho^k\right)-\widetilde N_{s,k}\right\|_{L^2(Y)\to L^2(Y)}
\lesssim h^{s+1}|\mu|^{-(3s+2-k)/2}.
\end{equation}
\end{lemma}

{\it Proof.} It follows from Lemma 2.6 that
$$\varrho^ka_{1,j}\in S^{-3j-2+k}(|\varrho|)\subset S^{-(3j+2-k)/2}(|\mu|),\quad j\ge s.$$
Therefore, by (\ref{eq:2.1}) we get
\begin{equation}\label{eq:2.26}
\left\|{\rm Op}_h\left(\varrho^ka_{1,j}\right)\right\|_{L^2(Y)\to L^2(Y)}\lesssim |\mu|^{-(3j+2-k)/2},\quad j\ge s,
\end{equation}
which clearly implies (\ref{eq:2.25}).
\eproof

\section{Parametrix construction in the glancing region} 

We will use the parametrix from the previous section to construct a parametrix for the boundary value
problem (\ref{eq:1.2}) with $f$ replaced by ${\rm Op}_h(\chi)f$. Then using our parametrix we will prove the following

\begin{Theorem} Let $|{\rm Im}\,z|\ge h^{2/3-\epsilon}$, $0<\epsilon\ll 1$. Then, for every integer $s\ge 0$ there is an operator
$${\cal A}_s(h,z)={\cal O}\left(h|{\rm Im}\,z|^{-1}\right):L^2(\partial X)\to L^2(\partial X)$$
independent of all $n_\ell$ with $\ell\ge s$ such that
\begin{equation}\label{eq:3.1}
\left\|N(h,z){\rm Op}_h(\chi)-\widetilde{{\rm Op}}_h(\rho\chi+c_sh^s\rho^{-s-1}zn_s\chi)-{\cal A}_s\right\|
\lesssim h^{s+1}|{\rm Im}\,z|^{-(3s+2)/2}
\end{equation}
where $c_s=0$, ${\cal A}_s=0$ if $s=0$, and $c_s=-i(-2i)^{-s-1}$ for $s\ge 1$. Furthermore, for every integer $s\ge 1$ there are operators
$${\cal B}_s^R(h,z), {\cal B}_s^L(h,z)={\cal O}\left(h^{-s}\right):L^2(\partial X)\to L^2(\partial X)$$
independent of all $n_\ell$ with $\ell\ge 1$, and operators
${\cal C}_s^R(h,z)$, ${\cal C}_s^L(h,z)$ independent of all $n_\ell$ with $\ell\ge s$ such that
\begin{equation}\label{eq:3.2}
\left\|N(h,z){\cal B}_s^R-{\cal C}_s^R-{\rm Op}_h(n_s\chi)\right\|
\lesssim h|{\rm Im}\,z|^{-(2s+1)/2},
\end{equation}
\begin{equation}\label{eq:3.3}
\left\|{\cal B}_s^LN(h,z)-{\cal C}_s^L-{\rm Op}_h(n_s\chi)\right\|
\lesssim h|{\rm Im}\,z|^{-(2s+1)/2}.
\end{equation}
\end{Theorem}

Clearly, Theorem 1.1 follows from (\ref{eq:3.1}) with $s=0$. Note that an analog of Theorem 3.1 has been proved in \cite{kn:V4} but for $|{\rm Im}\,z|\ge h^{1/2-\epsilon}$ and with worse bounds in the right-hand sides of (\ref{eq:3.1}), (\ref{eq:3.2}) and (\ref{eq:3.3}).

We begin the parametrix construction by writing the Laplace-Beltrami operator
in local coordinates near the boundary. Fix a point 
$x^0\in\Gamma$ and let ${\cal U}_0\subset\Gamma$ be a small open neighbourhood of $x^0$. Let $(x_1,x')$, $x_1>0$, $x'\in 
{\cal U}_0$, be the normal coordinates with respect to the metric ${\cal G}$. In these coordinates the Laplacian can be written as follows
$$\Delta_X=\partial_{x_1}^2+r(x,\partial_{x'})+q(x,\partial_x)$$
where $r(x,\xi')$ is homogeneous in $\xi'$ of order two, smooth in $x$ and strictly positive for all $\xi'\neq 0$, $q(x,\xi)=\langle q(x),\xi\rangle=q^\sharp(x)\xi_1+\langle q^\flat(x),\xi'\rangle$, $q^\sharp$ and $q^\flat$ being smooth functions. Clearly,
$r(0,x',\xi')=r_0(x',\xi')$ is the principal symbol of the Laplace-Beltrami operator on the boundary. Introduce the function
$$\varphi^\sharp(x_1,x')=-\frac{1}{2}\int_0^{x_1}q^\sharp(\sigma,x')d\sigma$$
and observe that
$$e^{-\varphi^\sharp}\Delta_X e^{\varphi^\sharp}=\partial_{x_1}^2+r(x,\partial_{x'})+q^\flat(x,\partial_{x'})+V^\sharp(x)$$
with a new function $q^\flat(x,\xi')=\langle q^\flat(x),\xi'\rangle$, $q^\flat$ and $V^\sharp$ being smooth functions.
We now introduce a new normal variable
$t=n_0(x')^{1/2}x_1$, $n_0(x'):=n(0,x')>0$, and we write the operator
$$P(h)=-h^2n_0^{-1}\Delta_X-1-i{\rm Im}\,z-zn_0^{-1}(n-n_0)$$
in the coordinates $(t,x')$ as follows:
$$P_\sharp(h):=e^{-\varphi^\sharp}P(h)e^{\varphi^\sharp}$$ $$={\cal D}_t^2+n_0^{-1}r_0(x',{\cal D}_{x'})-1-i{\rm Im}\,z-zn_0^{-1}\left(n(tn_0^{-1/2},x')-n_0\right)$$ 
$$+n_0^{-1}\left(r(tn_0^{-1/2},x',{\cal D}_{x'})-r(0,x',{\cal D}_{x'})\right)
-ihn_0^{-1}q^\flat(tn_0^{-1/2},x',{\cal D}_{x'})-h^2n_0^{-1}V^\sharp(tn_0^{-1/2},x').$$ 
If $u= e^{\varphi^\sharp}v$, we have $u|_{t=0}=v|_{t=0}$ and 
\begin{equation}\label{eq:3.4}
{\cal D}_{\nu}u|_{\partial X}={\cal D}_{x_1}u|_{x_1=0}={\cal D}_{x_1}v|_{x_1=0}+{\cal D}_{x_1}\varphi^\sharp|_{x_1=0} v|_{x_1=0}$$
$$=n_0^{1/2}{\cal D}_{t}v|_{t=0}+\frac{ih}{2}q^\sharp(0,x')v|_{t=0}.
\end{equation}
If we denote 
$$n_k(x')=\partial_\nu^kn|_{\partial X}=\partial_{x_1}^kn(x_1,x')|_{x_1=0},$$ $$  
V^\sharp_k(x')=\partial_\nu^kV^\sharp|_{\partial X}=\partial_{x_1}^kV^\sharp(x_1,x')|_{x_1=0},$$
$$r_k(x',\xi')=\partial_\nu^kr(x,\xi')|_{\partial X}=\partial_{x_1}^kr(x_1,x',\xi')|_{x_1=0},$$ $$ 
q^\flat_k(x',\xi')=\partial_\nu^kq^\flat(x,\xi')|_{\partial X}=\partial_{x_1}^kq^\flat(x_1,x',\xi')|_{x_1=0},$$
we have the formal expansions
$$n_0^{-1}\left(n(tn_0^{-1/2},x')-n_0\right)\sim\sum_{k=1}^\infty\frac{t^k}{k!}n_kn_0^{-(k+2)/2},$$
$$n_0^{-1}V^\sharp(tn_0^{-1/2},x')\sim\sum_{k=0}^\infty\frac{t^k}{k!}V^\sharp_kn_0^{-(k+2)/2},$$
$$n_0^{-1}\left(r(tn_0^{-1/2},x',\xi')-r(0,x',\xi')\right)\sim\sum_{k=1}^\infty\frac{t^k}{k!}r_kn_0^{-(k+2)/2},$$
$$n_0^{-1}q^\flat(tn_0^{-1/2},x',\xi')\sim\sum_{k=0}^\infty\frac{t^k}{k!}q^\flat_kn_0^{-(k+2)/2}.$$
 It is clear that it suffices to build the parametrix microlocally near $\Sigma$, that is, with
$\chi$ supported in a small neighbourhood, ${\cal W}_0\subset T^*\partial X$, of a point $\zeta^0=(x^0,\xi^0)\in\Sigma$. 
Then the global parametrix is obtained by making a suitable partition of the unity on $\Sigma$ as explained in the introduction. We may suppose that 
$\pi{\cal W}_0\subset{\cal U}_0$, where $\pi:T^*\partial X\to \partial X$ denotes the projection $\pi(x',\xi')=x'$.
Let $\psi\in C_0^\infty({\cal W}_0)$ be such that $\psi=1$ on supp$\chi$. The standard calculas of $h-\Psi$DOs
(e.g see Section 7 of \cite{kn:DS}) yield that to the function $\chi$ one can associate a linear map
$\vartheta:C^\infty(T^*\partial X)\to C_0^\infty({\cal W}_0)$ so that, if $a$ is a symbol independent of $h$, then 
$${\rm Op}_h(a){\rm Op}_h(\chi)-{\rm Op}_h(\chi){\rm Op}_h(a)={\rm Op}_h(\vartheta(a))+{\cal O}(h^\infty)$$
where 
$$\vartheta(a)\sim\sum_{j=1}^\infty h^j\vartheta_j(a)$$
with functions $\vartheta_j(a)$ independent of $h$ and supported in an arbitrary neighbourhood of supp$\,\chi$.
Clearly, we can rewrite the above identity in the form
$${\rm Op}_h(a){\rm Op}_h(\chi)={\rm Op}_h(\chi){\rm Op}_h(\psi a)+{\rm Op}_h(\vartheta(a))+{\cal O}(h^\infty).$$
 Using this we can write
\begin{equation}\label{eq:3.5}
P_\sharp{\rm Op}_h(\chi)={\rm Op}_h(\chi)\widetilde P_\sharp +{\cal O}(h^\infty)
\end{equation}
where
$$\widetilde P_\sharp={\cal D}_t^2+(\psi n_0^{-1}r_0)(x',{\cal D}_{x'})-1-i{\rm Im}\,z+p(t,x',{\cal D}_{x'};h,z)$$
with a function $p(t,x',\xi';h,z)={\cal O}(t+h)$ having the formal expansion
$$p\sim\sum_{k=0}^\infty\sum_{j=0}^\infty t^kh^jp_{k,j}$$
where the functions $p_{k,j}$ are independent of $h$ and $t$. More precisely, we have the formulas
$$p_{0,0}=0,\quad p_{k,0}=\frac{1}{k!}\psi (r_k-zn_k)n_0^{-(k+2)/2},$$
for $k\ge 1$,
$$p_{0,1}=-i\psi q_0^\flat n_0^{-1}+\vartheta_1(r_0/n_0),$$
$$p_{k,1}=-\frac{i}{k!}\psi q^\flat_kn_0^{-(k+2)/2}+
\frac{1}{k!}\vartheta_1\left((r_k-zn_k)n_0^{-(k+2)/2}\right),$$
for $k\ge 1$,
$$p_{0,2}=-\psi V_0^\sharp n_0^{-1}-i\vartheta_1(q_0^\flat/n_0)+\vartheta_2(r_0/n_0),$$
$$p_{k,2}=-\frac{1}{k!}\psi V^\sharp_kn_0^{-(k+2)/2}-\frac{i}{k!}\vartheta_1\left(q^\flat_kn_0^{-(k+2)/2}\right)+
\frac{1}{k!}\vartheta_2\left((r_k-zn_k)n_0^{-(k+2)/2}\right),$$
for $k\ge 1$,
$$p_{0,j}=-\vartheta_{j-2}\left(V_0^\sharp/n_0\right)-i\vartheta_{j-1}(q_0^\flat/n_0)+\vartheta_j(r_0/n_0),$$
for $j\ge 3$,
$$p_{k,j}=-\frac{1}{k!}\vartheta_{j-2}\left(V^\sharp_kn_0^{-(k+2)/2}\right)-\frac{i}{k!}\vartheta_{j-1}\left(q^\flat_kn_0^{-(k+2)/2}\right)+
\frac{1}{k!}\vartheta_j\left((r_k-zn_k)n_0^{-(k+2)/2}\right),$$
for $k\ge 1, j\ge 3$. We will now transform the operator $\widetilde P_\sharp$ into the normal form studied in the previous section by using
$h-$FIOs acting on the tangent variable $x'$ only and independent of the normal variable $t$. 
Roughly speaking, we have to make a suitable symplectic change of the variables $(x',\xi')$ so that in the new coordinates our operator has
a simpler form. Indeed, there exist an open,
bounded domain $Y\subset{\bf R}^{d-1}$ and a symplectomorphism $\kappa:{\cal W}_0\to T^*Y$ such that if $(y,\eta):=\kappa(x',\xi')$, then
$$\eta_1=\left(n_0^{-1}r_0\right)(x',\xi')-1$$
and $\kappa(\Sigma\cap{\cal W}_0)=\{\eta_1=0\}$. Let $U={\cal O}(1):L^2(\pi{\cal W}_0)\to L^2(Y)$ be an elliptic, zero-order $h-$FIO associated to
$\kappa$. Then the inverse $U^{-1}={\cal O}(1):L^2(Y)\to L^2(\pi{\cal W}_0)$ is an $h-$FIO associated to the inverse symplectomorphism 
$\kappa^{-1}$. It is well-known (e.g. see Section 11 of \cite{kn:Z}) that one can associate to $\kappa$ a linear map $\omega:C_0^\infty({\cal W}_0)\to C_0^\infty(T^*Y)$
so that, if $a$ is a symbol independent of $h$, then 
$$U{\rm Op}_h(a)U^{-1}={\rm Op}_h(\omega(a))+{\cal O}(h^\infty)$$
where 
$$\omega(a)\sim\sum_{j=0}^\infty h^j\omega_j(a)$$
with functions $\omega_j(a)$ independent of $h$, $\omega_0(a)=a\circ\kappa$. 
Moreover, each $\omega_j(a)$ is a linear combination of functions of the form 
$\left(\partial_{x'}^\alpha\partial_{\xi'}^\beta a\right)\circ\kappa$.
Therefore $\omega_j(a)$ is 
supported in an arbitrary neighbourhood of the set $\kappa({\rm supp}\,a)$. 
 Observe now that
$$(\psi n_0^{-1}r_0)\circ\kappa(y,\eta)=(\eta_1+1)\psi\circ\kappa(y,\eta)$$
and $\psi\circ\kappa=1$ on supp$\,\omega(\chi)$. Therefore, taking into account that ${\cal D}_{y_1}={\rm Op}_h(\eta_1)$,
we have modulo ${\cal O}(h^\infty)$,
$${\rm Op}_h(\chi)\left((\psi n_0^{-1}r_0)(x',{\cal D}_{x'})-1\right)U^{-1}$$
$$={\rm Op}_h(\chi)U^{-1}
\left({\rm Op}_h\left(\omega(\psi n_0^{-1}r_0)\right)-1\right)$$
$$={\rm Op}_h(\chi)U^{-1}
{\rm Op}_h\left(\omega(\psi n_0^{-1}r_0)-\omega_0(\psi n_0^{-1}r_0)\right)$$
$$+{\rm Op}_h(\chi)U^{-1}
\left({\rm Op}_h\left((\eta_1+1)\psi\circ\kappa\right)-1\right)$$
$$={\rm Op}_h(\chi)U^{-1}
{\rm Op}_h\left(\omega(\psi n_0^{-1}r_0)-\omega_0(\psi n_0^{-1}r_0)\right)$$
$$+{\rm Op}_h(\chi)U^{-1}
\left({\cal D}_{y_1}+{\rm Op}_h\left((\eta_1+1)(\psi\circ\kappa-1)\right)\right)$$
$$={\rm Op}_h(\chi)U^{-1}\left({\cal D}_{y_1}+
{\rm Op}_h\left(\omega(\psi n_0^{-1}r_0)-\omega_0(\psi n_0^{-1}r_0)\right)\right)$$
$$+U^{-1}{\rm Op}_h(\omega(\chi))
{\rm Op}_h\left((\eta_1+1)(\psi\circ\kappa-1)\right)$$
$$={\rm Op}_h(\chi)U^{-1}\left({\cal D}_{y_1}+
{\rm Op}_h\left(\omega(\psi n_0^{-1}r_0)-\omega_0(\psi n_0^{-1}r_0)\right)\right).$$
Thus by (\ref{eq:3.5}) we get
\begin{equation}\label{eq:3.6}
P_\sharp{\rm Op}_h(\chi)U^{-1}={\rm Op}_h(\chi)U^{-1}P_0+{\cal O}(h^\infty)
\end{equation}
where $P_0$ is the operator from the previous section with $\mu={\rm Im}\,z$ and
$$m=\omega\left(p+\psi n_0^{-1}r_0\right)-\omega_0\left(\psi n_0^{-1}r_0\right).$$
It is easy to see that $m$ has the formal expansion
$$m\sim\sum_{k=0}^\infty\sum_{j=0}^\infty t^kh^jm_{k,j}$$
where $m_{0,0}=0$ and
$$m_{0,j}=\omega_j\left(\psi n_0^{-1}r_0\right)+\sum_{\ell=0}^j\omega_{j-\ell}\left(p_{0,\ell}\right)$$
for $j\ge 1$,
$$m_{k,j}=\sum_{\ell=0}^j\omega_{j-\ell}\left(p_{k,\ell}\right)$$
for $k\ge 1$, $j\ge 0$. In particular, we have
\begin{equation}\label{eq:3.7}
m_{k,0}=\frac{1}{k!}\left((r_k-zn_k)n_0^{-(k+2)/2}\psi\right)\circ\kappa,\quad k\ge 1.
\end{equation}
Let $k\ge 0$ be an arbitrary integer and let $\widetilde u$ be the parametrix constructed in the previous section with 
$$\widetilde u|_{t=0}={\rm Op}_h(\phi\varrho^k)U{\rm Op}_h(\psi)f$$
where $\phi=\phi(\eta_1)$ is the same function as in (\ref{eq:2.3}).  
Taking supp$\,\psi$ small enough we can arrange that $\phi=1$ on supp$\,\omega(\psi)$.
If $\widetilde N$ is the operator from the previous section, then
$${\cal D}_t\widetilde u|_{t=0}=\widetilde N{\rm Op}_h(\phi_1\varrho^k)U{\rm Op}_h(\psi)f$$
where $\phi_1=\phi_1(\eta_1)\in C_0^\infty$ is such that $\phi_1=1$ on supp$\,\phi$. 
Set
$$v={\rm Op}_h(\chi)U^{-1}\widetilde u.$$ 
 We have
$$v|_{t=0}={\rm Op}_h(\chi)U^{-1}{\rm Op}_h(\phi\varrho^k)U{\rm Op}_h(\psi)f$$
and 
$${\cal D}_tv|_{t=0}={\rm Op}_h(\chi)U^{-1}\widetilde N{\rm Op}_h(\phi_1\varrho^k)U{\rm Op}_h(\psi)f.$$
If we set $w= e^{\varphi^\sharp}v$, then
\begin{equation}\label{eq:3.8}
w|_{\partial X}={\rm Op}_h(\chi)U^{-1}{\rm Op}_h(\phi\varrho^k)U{\rm Op}_h(\psi)f
\end{equation}
and, in view of (\ref{eq:3.4}),
\begin{equation}\label{eq:3.9}
{\cal D}_{\nu}w|_{\partial X}={\cal T}^{(k)}f
\end{equation}
where 
$${\cal T}^{(k)}=n_0^{1/2}{\rm Op}_h(\chi)U^{-1}\widetilde N{\rm Op}_h(\phi_1\varrho^k)U{\rm Op}_h(\psi)
+\frac{ih}{2}q^\sharp(0,x'){\rm Op}_h(\chi)U^{-1}{\rm Op}_h(\phi\varrho^k)U{\rm Op}_h(\psi)$$
$$=\sum_{j=0}^{M+1}h^j{\cal T}_j^{(k)}$$
where
$${\cal T}_0^{(k)}=n_0^{1/2}{\rm Op}_h(\chi)U^{-1}{\rm Op}_h(\phi\varrho^{k+1})U{\rm Op}_h(\psi),$$
$${\cal T}_1^{(k)}=\frac{i}{2}q^\sharp(0,x'){\rm Op}_h(\chi)U^{-1}{\rm Op}_h(\phi\varrho^k)U{\rm Op}_h(\psi)-in_0^{1/2}{\rm Op}_h(\chi)U^{-1}{\rm Op}_h(\varrho^ka_{1,0})U{\rm Op}_h(\psi),$$
$${\cal T}_j^{(k)}=-in_0^{1/2}{\rm Op}_h(\chi)U^{-1}{\rm Op}_h(\varrho^ka_{1,j-1})U{\rm Op}_h(\psi),$$
for $j\ge 2$. Set
$${\cal H}_j^{(k)}=n_0^{1/2}{\rm Op}_h(\chi)U^{-1}{\rm Op}_h\left(\varrho^{k-j-1}\left(n_jn_0^{-(j+2)/2}\psi\right)\circ\kappa\right)U{\rm Op}_h(\psi),$$
$${\cal K}_j^{(k)}=n_jn_0^{-(j+1)/2}{\rm Op}_h(\chi)U^{-1}{\rm Op}_h\left(\phi\varrho^{k-j-1}\right)U{\rm Op}_h(\psi).$$

\begin{lemma} For $0\le k\le j+1$ we have the bound
\begin{equation}\label{eq:3.10}
\left\|{\cal H}_j^{(k)}-{\cal K}_j^{(k)}\right\|\lesssim
h|{\rm Im}\,z|^{-(j+1-k)/2}.
\end{equation}
\end{lemma}

{\it Proof.} Clearly
$$\varrho^{k-j-1}\left(n_jn_0^{-(j+2)/2}\psi\right)\circ\kappa,\,\,\phi\varrho^{k-j-1}\in S^{-(j+1-k)}(|\varrho|)
\subset S^{-(j+1-k)/2}(|\mu|)$$
and hence
$${\rm Op}_h\left(\varrho^{k-j-1}\left(n_jn_0^{-(j+2)/2}\psi\right)\circ\kappa\right),\,\,
{\rm Op}_h\left(\phi\varrho^{k-j-1}\right)={\cal O}\left(|{\rm Im}\,z|^{-(j+1-k)/2}\right)$$
in the $L^2(Y)\to L^2(Y)$ norm. On the other hand, we have, mod ${\cal O}(h)$,
$$n_jn_0^{-(j+1)/2}{\rm Op}_h(\chi)U^{-1}=n_0^{1/2}{\rm Op}_h(\chi){\rm Op}_h\left(n_jn_0^{-(j+2)/2}\psi\right)U^{-1}$$
$$=n_0^{1/2}{\rm Op}_h(\chi)U^{-1}\left({\rm Op}_h\left(\left(n_jn_0^{-(j+2)/2}\psi\right)\circ\kappa\right)+{\cal O(h)}\right).$$
Therefore, modulo an operator of norm ${\cal O}\left(h|{\rm Im}\,z|^{-(j+1-k)/2}\right)$, we obtain
$${\cal H}_j^{(k)}-{\cal K}_j^{(k)}=n_0^{1/2}{\rm Op}_h(\chi)U^{-1}\left({\rm Op}_h\left(\varrho^{k-j-1}\left( n_jn_0^{-(j+2)/2}\psi\right)\circ\kappa\right)\right.$$ $$-\left.
{\rm Op}_h\left(\left(n_jn_0^{-(j+2)/2}\psi\right)\circ\kappa\right){\rm Op}_h\left(\phi\varrho^{k-j-1}\right)\right)U{\rm Op}_h(\psi).$$
Since $\phi=1$ on supp$\,\psi\circ\kappa$, it is easy to see that the bound (\ref{eq:3.10}) follows from Proposition 2.1.
\eproof

\begin{lemma} For every $j\ge 1$ the operator
${\cal T}_j^{(k)}+i(-2i)^{-j-1}z{\cal H}_j^{(k)}$
is independent of all $n_\ell$ with $\ell\ge j$.
\end{lemma}

{\it Proof.} By Lemma 2.6 the function
$$a_{1,j-1}+\frac{j!\phi m_{j,0}}{(-2i\varrho)^{j+1}}$$
depends only on $m_{\nu,0}$ with $\nu\le j-1$ and $m_{\nu,\ell}$ with $1\le\ell\le j$, $\nu\le j-\ell\le j-1$. On the other hand, it is clear from
the above formulas that $m_{\nu,\ell}$ with $\nu\le j-1$ depends only on $n_\ell$ with $\ell\le j-1$. Thus by (\ref{eq:3.7}) we conclude that
$$a_{1,j-1}-z(-2i\varrho)^{-j-1}\left(n_jn_0^{-(j+2)/2}\psi\right)\circ\kappa$$
is independent of all $n_\ell$ with $\ell\ge j$, which clearly implies the lemma.
\eproof

Proposition 2.8 and Lemma 2.9 imply the following

\begin{prop} Given any integers $s,k\ge 0$ such that $k\le 3s+2$, we have the estimate
\begin{equation}\label{eq:3.11}
\left\|N(h,\mu){\rm Op}_h(\chi)U^{-1}{\rm Op}_h(\phi\varrho^k)U{\rm Op}_h(\psi)-\sum_{j=0}^{s}h^j{\cal T}_j^{(k)}\right\|
\lesssim h^{s+1}|{\rm Im}\,z|^{-(3s+2-k)/2}.
\end{equation}
For $j\ge 1$, $0\le k\le 3j-1$ we also have
\begin{equation}\label{eq:3.12}
\left\|{\cal T}_j^{(k)}\right\|\lesssim |{\rm Im}\,z|^{-(3j-1-k)/2}.
\end{equation}
\end{prop}

{\it Proof.} Clearly (\ref{eq:3.12}) follows from (\ref{eq:2.26}) used with $j$ replaced by $j-1$.
To prove (\ref{eq:3.11}) we will make use of the coercivity of the Dirichlet realization, $G_D$, of the operator
$n^{-1}\Delta_X$ in the Hilbert space $L^2(X;ndx)$. We have
$$\|u\|_{H_h^2(X)}\lesssim \|h^2G_Du\|_{L^2(X)}+\|u\|_{L^2(X)},\quad\forall u\in D(G_D),$$
where $H_h^2(X)$ denotes the Sobolev space equipped with the semi-classical norm. This together with the semi-classical
version of the trace theorem imply
\begin{equation}\label{eq:3.13}
\left\|\gamma{\cal D}_\nu\left(h^2G_D+z\right)^{-1}\right\|_{L^2(X)\to L^2(\partial X)}\lesssim h^{-1/2}
\left\|\left(h^2G_D+z\right)^{-1}\right\|_{L^2(X)\to H_h^2(X)}$$
$$\lesssim h^{-1/2}+h^{-1/2}\left\|\left(h^2G_D+z\right)^{-1}\right\|_{L^2(X)\to L^2(X)}\lesssim h^{-1/2}|{\rm Im}\,z|^{-1}\lesssim h^{-7/6}
\end{equation}
where $\gamma$ denotes the restriction on $\partial X$. Let $u$ be the solution to equation (\ref{eq:1.2}) with boundary condition
$$u|_{\partial X}={\rm Op}_h(\chi)U^{-1}{\rm Op}_h(\phi\varrho^k)U{\rm Op}_h(\psi)f.$$
If $v$, $w$ and $\widetilde u$ are the functions introduced above, we have in view of (\ref{eq:3.6}) and (\ref{eq:3.8}), that
$(u-w)|_{\partial X}=0$ and
$$\left(h^2n^{-1}\Delta_X+z\right)(u-w)=P(h)w=P_\sharp(h)v={\rm Op}_h(\chi)U^{-1}P_0\widetilde u+{\cal O}(h^\infty)\widetilde u.$$
Hence
$$u-w=\left(h^2G_D+z\right)^{-1}\left({\rm Op}_h(\chi)U^{-1}P_0+{\cal O}(h^\infty)\right)\widetilde u$$
which together with (\ref{eq:3.9}) imply the identity
\begin{equation}\label{eq:3.14}
N{\rm Op}_h(\chi)U^{-1}{\rm Op}_h(\phi\varrho^k)U{\rm Op}_h(\psi)f-{\cal T}^{(k)}f$$ $$
=\gamma{\cal D}_\nu\left(h^2G_D+z\right)^{-1}\left({\rm Op}_h(\chi)U^{-1}P_0+{\cal O}(h^\infty)\right)\widetilde u.
\end{equation}
By Proposition 2.8, (\ref{eq:3.13}) and (\ref{eq:3.14}) we get
$$
\left\|N{\rm Op}_h(\chi)U^{-1}{\rm Op}_h(\phi\varrho^k)U{\rm Op}_h(\psi)f-{\cal T}^{(k)}f\right\|_{L^2(\partial X)}$$ 
\begin{equation}\label{eq:3.15}
\lesssim h^{-7/6}\left\|P_0\widetilde u\right\|_{L^2((0,1)\times Y)}+{\cal O}(h^\infty)\left\|\widetilde u\right\|_{L^2((0,1)\times Y)}
\end{equation}
 $$\lesssim h^{\epsilon M/2-7/6}\left\|{\rm Op}_h(\phi\varrho^k)U{\rm Op}_h(\psi)f\right\|_{L^2(Y)}
\lesssim h^{\epsilon M/2-7/6}\left\|f\right\|_{L^2(\partial X)}.$$
On the other hand, by Lemma 2.9 we have
\begin{equation}\label{eq:3.16}
\left\|{\cal T}^{(k)}-\sum_{j=0}^{s}h^j{\cal T}_j^{(k)}\right\|\lesssim h^{s+1}|{\rm Im}\,z|^{-(3s+2-k)/2}.
\end{equation}
Take now $M$ so that $\epsilon M/2-7/6>s+1$. 
Clearly, the bound (\ref{eq:3.11}) follows from (\ref{eq:3.15}) and (\ref{eq:3.16}).
\eproof

It follows from Lemma 3.3 that for $s\ge 1$ the operator
$${\cal P}_s^{(k)}=\sum_{j=0}^{s}h^j{\cal T}_j^{(k)}+i(-2i)^{-s-1}zh^s{\cal H}_s^{(k)}$$
is independent of all $n_\ell$ with $\ell\ge s$. Put ${\cal P}_0^{(0)}={\cal T}_0^{(0)}$.
By (\ref{eq:3.10}) and (\ref{eq:3.11}) we have for $s=k=0$ and $s\ge 1$, $0\le k\le s+1$,
\begin{equation}\label{eq:3.17}
\left\|N{\rm Op}_h(\chi)U^{-1}{\rm Op}_h(\phi\varrho^k)U{\rm Op}_h(\psi)-{\cal P}_s^{(k)}-zc_sh^s{\cal K}_s^{(k)}\right\|
\lesssim h^{s+1}|{\rm Im}\,z|^{-(3s+2-k)/2}
\end{equation}
where $c_s$ is as in Theorem 3.1. 
Since
$${\rm Op}_h(\chi)U^{-1}{\rm Op}_h(\phi)U{\rm Op}_h(\psi)={\rm Op}_h(\chi)+{\cal O}(h^\infty),$$
it is easy to see that (\ref{eq:3.1}) follows from (\ref{eq:3.17}) with $k=0$ and ${\cal A}_s={\cal P}_s^{(0)}-{\cal T}_0^{(0)}$.
Since
$$n_sn_0^{-(s+1)/2}{\rm Op}_h(\chi)U^{-1}{\rm Op}_h(\phi)U{\rm Op}_h(\psi)$$ $$
=n_0^{-(s+1)/2}{\rm Op}_h(n_s\chi)+{\cal O}(h)
={\rm Op}_h(n_s\chi)n_0^{-(s+1)/2}+{\cal O}(h),$$
the bound (\ref{eq:3.2}) follows from (\ref{eq:3.17}) with $k=s+1$ and 
$${\cal B}_s^R=\left(zc_sh^s\right)^{-1}{\rm Op}_h(\chi)U^{-1}{\rm Op}_h(\phi\varrho^{s+1})U{\rm Op}_h(\psi)n_0^{(s+1)/2},$$
 $${\cal C}_s^R=\left(zc_sh^s\right)^{-1}{\cal P}_s^{(s+1)}n_0^{(s+1)/2}.$$
Since
$$N(h,z)=N(h,\overline z)^*$$
and
$${\rm Op}_h(n_s\chi)=\left({\rm Op}_h(n_s\chi)\right)^*+{\cal O}(h),$$
the bound (\ref{eq:3.3}) follows from (\ref{eq:3.2}) with 
$${\cal B}_s^L(h,z)={\cal B}_s^R(h,\overline z)^*,\quad {\cal C}_s^L(h,z)={\cal C}_s^R(h,\overline z)^*.$$ 
The analysis of the operator $N{\rm Op}_h(1-\chi)$ is much easier and, as mentioned in the introduction, can be done for
$|{\rm Im}\,z|\ge h^{1-\epsilon}$. Indeed, in this case the parametrix construction in Section 3 of \cite{kn:V4} gives full
expansions similar to those in Theorem 3.1 but with better bounds in the right-hand sides (with $|{\rm Im}\,z|$ replaced by $1$).
Thus, combining the results of \cite{kn:V4} with Theorem 3.1 we get the following

\begin{Theorem} Let $|{\rm Im}\,z|\ge h^{2/3-\epsilon}$, $0<\epsilon\ll 1$. Then, for every integer $s\ge 0$ there is an operator
$${\cal A}_s(h,z)={\cal O}\left(h|{\rm Im}\,z|^{-1}\right):L^2(\partial X)\to L^2(\partial X)$$
independent of all $n_\ell$ with $\ell\ge s$ such that
\begin{equation}\label{eq:3.18}
\left\|N(h,z)-\widetilde{{\rm Op}}_h(\rho+c_sh^s\rho^{-s-1}zn_s)-{\cal A}_s\right\|
\lesssim h^{s+1}|{\rm Im}\,z|^{-(3s+2)/2}
\end{equation}
where $c_s=0$, ${\cal A}_s=0$ if $s=0$, and $c_s=-i(-2i)^{-s-1}$ for $s\ge 1$. Furthermore, for every integer $s\ge 1$ there are operators
${\cal B}_s^R(h,z), {\cal B}_s^L(h,z)$
independent of all $n_\ell$ with $\ell\ge 1$, and operators
${\cal C}_s^R(h,z)$, ${\cal C}_s^L(h,z)$ independent of all $n_\ell$ with $\ell\ge s$ such that
\begin{equation}\label{eq:3.19}
\left\|N(h,z){\cal B}_s^R-{\cal C}_s^R-n_sI\right\|
\lesssim h|{\rm Im}\,z|^{-(2s+1)/2},
\end{equation}
\begin{equation}\label{eq:3.20}
\left\|{\cal B}_s^LN(h,z)-{\cal C}_s^L-n_sI\right\|
\lesssim h|{\rm Im}\,z|^{-(2s+1)/2},
\end{equation}
where $I$ denotes the identity.
\end{Theorem}

\section{Applications to the transmission eigenvalues}

Let $\Omega\subset{\bf R}^d$, $d\ge 2$, be a bounded, connected domain with a $C^\infty$ smooth boundary $\Gamma=\partial\Omega$. 
A complex number $\lambda\neq 0$, ${\rm Re}\,\lambda\ge 0$, will be said to be a transmission eigenvalue if the following problem has a non-trivial solution:
\begin{equation}\label{eq:4.1}
\left\{
\begin{array}{lll}
\left(\Delta+\lambda^2 n_1(x)\right)u_1=0 &\mbox{in} &\Omega,\\
\left(\Delta+\lambda^2 n_2(x)\right)u_2=0 &\mbox{in} &\Omega,\\
u_1=u_2,\,\,\, \partial_\nu u_1=\partial_\nu u_2& \mbox{on}& \Gamma,
\end{array}
\right.
\end{equation}
where $\nu$ denotes the Euclidean unit inner normal to $\Gamma$, $n_j\in C^\infty(\overline\Omega)$, 
$j=1,2$ are strictly positive real-valued functions. We have the following

\begin{Theorem} Suppose that there is an integer $j\ge 1$ such that
\begin{equation}\label{eq:4.2}
\partial_\nu^sn_1(x)\equiv \partial_\nu^sn_2(x)\quad\mbox{on}\quad\Gamma,\quad 0\le s\le j-1,
\end{equation}
\begin{equation}\label{eq:4.3}
\partial_\nu^jn_1(x)\neq \partial_\nu^jn_2(x)\quad\mbox{on}\quad\Gamma.
\end{equation}
Then there exists a constant $C>0$ such that there are no transmission eigenvalues in the region 
\begin{equation}\label{eq:4.4}
\left\{\lambda\in{\bf C}:{\rm Re}\,\lambda\ge 0,\,\,|{\rm Im}\,\lambda|\ge C\left({\rm Re}\,\lambda+1\right)^{1-k_j
}\right\},
\end{equation}
where $k_1=2/3-\epsilon$, $\forall 0<\epsilon\ll 1$, and $k_j=2(2j+1)^{-1}$ if $j\ge 2$.
\end{Theorem}

Note that this theorem has been proved in \cite{kn:V4} with $k_j=2(3j+2)^{-1}$, $j\ge 1$. Here we get
a larger eigenvalue-free region. Previously, smaller eigenvalue-free regions were obtained in \cite{kn:LV}. 
Theorem 4.1 is an immediate consequence of the estimate (\ref{eq:3.20}) above. The proof goes in precisely the same way
as in Section 5 of \cite{kn:V4} using (\ref{eq:3.20}) instead of Theorem 4.1 of \cite{kn:V4}, and therefore we omit the details here.

In the non-degenerate isotropic case when
\begin{equation}\label{eq:4.5}
n_1(x)\neq n_2(x)\quad\mbox{on}\quad\Gamma
\end{equation}
it has been proved in \cite{kn:V3} that there are no transmission eigenvalues in a much larger region of the form
\begin{equation}\label{eq:4.6}
\left\{\lambda\in{\bf C}:{\rm Re}\,\lambda\ge 0,\,\,|{\rm Im}\,\lambda|\ge C\right\}
\end{equation}
for some constant $C>0$, which is in fact the optimal eigenvalue-free region. 
Note that it follows from \cite{kn:CLM} (see Theorem 4.2) that under the condition (\ref{eq:4.2}) 
the eigenvalue-free region (\ref{eq:4.6}) is no longer valid. Note also that parabolic eigenvalue-free regions
imply Weyl asymptotics for the counting function of the transmission eigenvalues with remainder term depending on the size
of the eigenvalue-free region (see \cite{kn:PV}). Roughly, the larger the eigenvalue-free region is, the smaller the 
remainder term is.

\end{document}